\documentclass[10pt]{article}

\usepackage{a4wide}
\usepackage{amssymb}
\usepackage{amsfonts}
\usepackage{amsmath}
\input xy
\xyoption{arrow} \xyoption{matrix}

\date{}

\newtheorem{proposition}{Proposition}[section]
\newtheorem{theorem}[proposition]{Theorem}
\newtheorem{lemma}[proposition]{Lemma}

\newtheorem{corollary}[proposition]{Corollary}

\def\der{\partial }

\def\nFM0{{\nu }_{F,M_0}}
\def\nFN0{{\nu }_{F,N_0}}
\def\nGN0{{\nu }_{G,N_0}}

\def\N0{ {\bf N}_0 }

\def\t{\otimes}
\def\g{\gamma}

\def\ra{\rightarrow}

\def\lra{\leftrightarrow}
\def\Xpm{X^{\pm }}

\def\s{\sigma}
\def\Z{{\bf Z }}

\def\l1{{\lambda}_1}

\def\a{\alpha}
\def\a0{ {\alpha }_0}
\def\a1{ {\alpha }_1}

\def\l{\lambda}
\def\o{\omega}
\def\nFGM0{{\nu }_{F,G,M_0}}

\def\nFN0{{\nu}_{F,N_0}}

\def\sm{{\sigma}^m}

\def\sm1{{\sigma}^{-1}}

\def\smtp1{{\sigma}^{-t+1}}

\def\o{\omega }
\def\S1{S^{-1}}

\def\Xpm1{X^{\pm 1}_1}

\def\sPM1{{\sigma }^{\pm 1}}
\def\sMP1{{\sigma }^{\mp 1 }}

\def\d{\delta}

\def\di{{\rm d.ind}}

\def\L{\Lambda}

\def\CA{{\cal A}}

\def\CD{{\cal D}}

\def\Ytm1{Y^{t-1}}
\def\Yim1{Y^{i-1}}

\def\CM{{\cal M}}
\def\CN{{\cal N}}

\def\CG{{\cal G}}
\def\CH{{\cal H}}

\def\supp{{\rm supp }}

\def\ad{{\rm ad }}
\def\dim{{\rm dim }}

\def\ker{ {\rm ker } }

\def\SL2Z{ {\rm SL}_2({\bf Z}) }

\def\Gp1{ G^{1 , 1 } }
\def\P11{ P^{-1 , 1 } }
\def\Pp1{ P^{1 , 1 } }

\def\nCLsr{{}^\nu\kern-2pt {\cal L}^{\sigma , \rho  }}
\def\nP{{}^\nu \kern-2pt P}
\def\nL{{}^\nu\kern-2pt L}
\def\nLL{{}^\nu\kern-2pt \Lambda}
\def\nPsr{{}^\nu\kern-2pt P^{\sigma , \rho  }}
\def\nLsr{{}^\nu\kern-2pt L^{\sigma , \rho  }}
\def\nuCL{{}^\nu\kern-2pt  {\cal L}}
\def\nCLsr{{}^\nu\kern-2pt {\cal L}^{\sigma , \rho  }}
\def\nCL1m{{}^\nu\kern-2pt {\cal L}^{-1 , 1  }}
\def\x1nu{x^\frac{1}{\nu}}
\def\xm1nu{x^{-\frac{1}{\nu}}}

\def\CN{{\cal N}}
\def\ra{\rightarrow }

\def\CB{{\cal B}}

\def\CC{ {\cal C}}

\def\CH{ {\cal H}}

\def\nAM0{{\nu }_{{\cal A},M_0}}
\def\nAN0{{\nu }_{{\cal A},N_0}}

\def\ad{ {\rm ad }}

\def\bx{\overline{x}}
\def\by{\overline{y}}

\def\gn{\mathfrak{n}}

\def\di!{\frac{\der^i}{i!}}
\def\dik!{\frac{\der^k_i}{k!}}

\def\CC{{\cal C}}

\def\id{{\rm id}}

\def\N{\mathbb{N}}

\def\0{\overline{0}}
\def\1{\overline{1}}

\def\Ln1{\L_{n,\overline{1}}}

\def\a1{a_{\overline{1}}}

\def\S{\Sigma}

\def\vn1{\overrightarrow{n-1}}

\def\Sh{{\rm Sh}}

\def\im{{\rm im}}

\def\mA{\mathbb{A}}

\def\mJ{\mathbb{J}}
\def\mI{\mathbb{I}}

\def\ann{{\rm ann}}

\def\mT{\mathbb{T}}

\def\K1{{\rm K}_1}

\def\hmI1{\widehat{\mI_1}}
\def\tmI1{\widetilde{\mI_1}}
\def\tmJ1{\widetilde{\mJ_1}}
\def\hB1{\widehat{B_1}}
\def\hCB1{\widehat{\CB_1}}

\def\Z{\mathbb{Z}}

\def\mB{\mathbb{B}}
\def\by{\overline{y}}
\def\bz{\overline{z}}

\begin{document}

\author{V. V. \  Bavula %(BiqAlg-3gen.tex)
}

\title{Description of  bi-quadratic algebras on 3 generators with PBW basis}

\maketitle
\begin{abstract}

The aim of the paper is to give an explicit  description of bi-quadratic algebras on 3 generators with PBW basis. \\
%Also a classification of bi-quadratic algebras on 2 generators is given (there are precisely 5 non-isomorphic algebras). 


$${\bf Contents}$$
\begin{enumerate}
\item Introduction.
\item Classification of bi-quadratic algebras on 2 generators.
\item Bi-quadratic algebras on 3 generators.
  \item The bi-quadratic algebras with $q_1\neq 1$, $q_2=q_3=1$.
   \item The bi-quadratic algebras with $q_1\neq 1$, $q_2\neq 1$, $q_3=1$.
    \item The bi-quadratic algebras with $q_1\neq 1$, $q_2\neq 1$, $q_3\neq 1$.
     \item The bi-quadratic algebras with $q_1=q_2=q_3=1$.
 % \item 
% \item 
%\item 
\end{enumerate}

\end{abstract}

%%%%%%%%%%%%%%%%%% SECTION 1 %%%%%%%%%%%%%%%%%%%%%%%%

\section{Introduction}

The aim of the paper is to give an explicit  description of bi-quadratic algebras on 3 generators with PBW basis. First we give necessary definitions and then explain the results and main ideas. \\

{\bf The skew bi-quadratic algebras.} For a ring $D$ and a natural number $n\geq 2$, a family $M=(m_{ij})_{i>j}$ of elements $m_{ij}\in D$ (where $1\leq j <i\leq n$) is called a {\em lower triangular half-matrix} with coefficients in $D$. The set of all such matrices is denoted by $L_n(D)$. \\

{\it Definition.}  Let $D$ be a ring, $Z (D)$ be its centre, $\sigma=(\sigma_1,\ldots  ,\sigma_n)$  be an $n$-tuple of
commuting endomorphisms of $D$, $\d=(\d_1,\ldots  ,\d_n)$ be   an $n$-tuple of
 $\s$-endomorphisms of $D$ ($\d_i$ is a $\s_i$-derivation of $D$ for $i=1, \ldots , n$), 
$Q=(q_{ij})\in L_n(Z(D))$, $\mA = (a_{ij,k})$ where $a_{ij,k}\in D$, $1\leq j <i\leq n$ and $k=1, \ldots , n$ and $\mB = (b_{ij})\in L_n(D)$. 
 The {\bf skew bi-quadratic algebra} (SBQA) $A= D[x_1, \ldots , x_n;  \s , \d , Q, \mA , \mB ]$ is a ring generated by the ring $D$ and elements $x_1, \ldots , x_n$ subject to the defining relations: For all $d\in D$, 
%\marginpar{SBQA-1}
\begin{equation}\label{SBQA-1}
x_id= \s_i(d) x_i+\d_i(d) \;\; {\rm for }\;\; i=1,\ldots , n, 
\end{equation}

%\marginpar{SBQA-2}
\begin{equation}\label{SBQA-2}
x_ix_j-q_{ij}x_jx_i=\sum_{k=1}^na_{ij,k}x_k+b_{ij}\;\; {\rm for \; all}\;\; j<i.
\end{equation}

{\it Definition.} In the particular case, when $\s_i =\id_D$ and $\d_i=0$ for $i=1,\ldots , n$, the ring $A$ is called the {\bf  bi-quadratic algebra} (BQA) and is denoted by  $A= D[x_1, \ldots , x_n;   Q, \mA , \mB ]$.\\

We say that the algebra $A= D[x_1, \ldots , x_n;   Q, \mA , \mB ]$ has {\em PBW basis} if $A= \bigoplus_{\alpha \in \N^n} Dx^\alpha$ where $x^\alpha = x_1^{\alpha_1}\cdots x_n^{\alpha_n}$. 

On the set $W_n$ of all words in the alphabet $\{ x_1, \ldots , x_n\}$ ($W_n$ is a multiplicative monoid freely generated by the elements $ x_1, \ldots , x_n$) we consider the {\em degree-by-lexicographic ordering} where $x_1<\cdots <x_n$. In more detail, $x_{i_1}\cdots x_{i_s}<x_{j_1}\cdots x_{j_t}$ if either $s<t$ or $s=t$, $i_1=j_1, \ldots ,  i_k= j_k$ and $i_{k+1}<j_{k+1}$ for some $k$ such that $1\leq k <s$. 

Given a bi-quadratic algebra $A$ on $n\geq 3$ generators. For each triple $i,j,k\in  \{ 1, \ldots , n\}$ such that $i<j<k$, there exactly two different ways to simplify the product $x_kx_jx_i$ with respect to the degree-by-lexicographic ordering:
%\marginpar{kjiL}
\begin{equation}\label{kjiL}
x_kx_jx_i= q_{kj}q_{ki}q_{ji}x_ix_jx_k+\sum_{|\alpha |\leq 2} c_{k,j,i,\alpha}x^\alpha, 
\end{equation}
%\marginpar{kjiR}
\begin{equation}\label{kjiR}
x_kx_jx_i= q_{kj}q_{ki}q_{ji}x_ix_jx_k+\sum_{|\alpha |\leq 2} c'_{k,j,i,\alpha}x^\alpha, 
\end{equation}
where in the first (resp. second) equality we start to simplify the product with the relation $x_kx_j=q_{kj}x_jx_k+\cdots $ (resp., $x_jx_i=q_{ji}x_ix_j+\cdots$). 

Let $S_n$ be the symmetric group of order $n$ (the group of all bijections from a set that contains $n$ elements to itself).

\begin{theorem}\label{3Jul16}%\marginpar{3Jul16}
Let  $A= D[x_1, \ldots , x_n;   Q, \mA , \mB ]$ be a bi-quadratic algebra where $n\geq 3$. Then the defining relations (\ref{SBQA-2}) are consistent and $A= \bigoplus_{\alpha \in \N^n} Dx^\alpha$  iff the following equalities hold: For all triples $i,j,k \in \{ 1, \ldots , n\}$ such that $i<j<k$, 
%\marginpar{ckjia}
\begin{equation}\label{ckjia}
c_{k,j,i,\alpha}=c_{k,j,i,\alpha}'.
\end{equation}
If (\ref{ckjia}) holds then,  for all $\s \in S_n$,  $A= \bigoplus_{\alpha \in \N^n} Dx_\s^\alpha$ where $x_\s^\alpha = x_{\s (1)}^{\alpha_1}\cdots x_{\s (n)}^{\alpha_n}$ and $\alpha = (\alpha_1, \ldots , \alpha_n)$. 
\end{theorem}

{\bf The groups $G_n= \Sh_n\rtimes \mT^n$ and $G_n'=G_n\rtimes S_n$.} 
 The following groups of linear changes of the variables $x_1, \ldots , x_n$ preserve the structure of the defining relations of the bi-quadratic algebras $A= K[x_1, \ldots , x_n; Q, \mA , \mB]$:
 \begin{eqnarray*}
 \mT^n &=& \{ t_\l \, |,\, \l = (\l_1, \ldots , \l_n) \in K^{\times n }\} \simeq ( K^{\times n }, \cdot ), \;\;  t_\l \lra \l \; \text{where $t_\l (x_i) = \l_ix_i$ for $i=1, \ldots , n$},\\
\Sh_n &=& \{ a_\mu \, |,\, \mu = (\mu_1, \ldots , \mu_n) \in K^n \} \simeq ( K^n, + ), \;\;  s_\mu \lra \mu \; \text{where $s_\mu (x_i) = x_i+\mu_i$ for $i=1, \ldots , n$}.
\end{eqnarray*} 
The groups $\mT^n$ and $\Sh_n$ are called the {\em algebraic $n$-dimensional  torus} and the {\em shift group}, respectively. 
 Clearly, 
%\marginpar{tlsmt}
\begin{equation}\label{tlsmt}
t_\l s_\mu t_\l^{-1} (x_i) = x_i+\l_i^{-1}\mu_i\;\; {\rm for}\;\; i=1,.\ldots , n.
\end{equation} 
 So, $t_\l s_\mu t_\l^{-1}\in \Sh_n$, $\mT^n\cap \Sh_n=\{ e\}$ and the group which is generated by $\mT^n$ and $\Sh_n$ is a semi-direct product of groups $G_n:=  \Sh_n\rtimes\mT^n$. Under the action of the group $G_n$ the half-matrix  $Q$ does not change. 
 
 The symmetric group $S_n$ acts on the set $\{ x_1, \ldots , x_n\}$ by permuting the variables: For all $\s \in S_n$ and $x_i$, $\s (x_i) = x_{\s (i)}$. The group $G_n'$, which is generated by $G_n$ and $S_n$, is a semi-direct product $G_n'=G_n\rtimes S_n$. The group $G_n$ acts on the class of bi-quadratic algebras. The group $G_n'$ acts on the class of bi-quadratic algebras $A=D[x_1, \ldots , x_n; Q, \mA , \mB ]$ where all elements of $Q$ are {\em invertible}: For $\s \in G_n'$, 
 $$ {}^\s A = D[\s (x_1), \ldots ,\s (x_n);   {}^\s Q,  {}^\s \mA ,  {}^\s \mB ]$$
 where the triple ${}^\s Q,  {}^\s \mA ,  {}^\s \mB $ can be explicitly written. In particular, for an element $\s$ of $S_n$, 
 $$ {}^\s A = D[x_{\s (1)}, \ldots , x_{\s (n)};   {}^\s Q,  {}^\s \mA ,  {}^\s \mB ].$$
For each half-matrix $Q=(q_{ij})_{i>j}$ with invertible entries $q_{ij}$, there is its unique completion $Q^c$ to an $n\times n$ matrix with 1 on the  diagonal (i.e. $q_{ii}=1$ for all $i=1, \ldots ,n$) and $q_{ji}=q_{ij}^{-1}$ for all $i>j$. 

The symmetric group $S_n$ acts on the set of $n\times n$ matrices by  simultaneously permuting their rows and columns (for $\s \in S_n$ and a matrix $M=(m_{ij})$, ${}^\s M = ({}^\s m_{ij})$ where ${}^\s m_{ij} := m_{\s (i)\s(j)}$). The matrix  ${}^\s Q$ of the algebra ${}^\s A$ is the lower triangular part of the matrix ${}^\s (Q^c)$. Clearly, for all $\s \in G_n$, ${}^\s Q=Q$.

When we say `up to $G_n$  (resp., $G_n'$)' we mean `up to affine change of the canonical generators $x_1, \ldots , x_n$ of the bi-quadratic algebra $A$  induced by an element of the group $G_n$ (resp., $G_n'$)'.\\

{\bf The bi-quadratic algebra of Lie type.}\\

{\it Definition.} A bi-quadratic algebra $A$ is called a {\bf bi-quadratic algebra of Lie type} provided all $q_{ij}=1$ and is denoted by  $A= D[x_1, \ldots , x_n;   \mA , \mB ]$. \\

Theorem\ref{6Jul16} explains why such a  name is given for this class of bi-quadratic algebras. 
\begin{theorem}\label{6Jul16}%\marginpar{6Jul16}
Let  $A= D[x_1, \ldots , x_n; Q,   \mA , \mB ]$ be a bi-quadratic algebra such that  $D$ is a $K$-algebra over a field $K$ and all the elements $a_{ij,k}$ and $b_{ij}$ belong to $K$.  Then the algebra $A$ is a bi-quadratic algebra of Lie type iff $A\simeq D\t_K(U(\CG )/(Z-1))$ where $\CG$ is a Lie algebra of dimension $n+1$ with nontrivial centre $Z(\CG )$  and  $Z\in \Z(\CG )\backslash \{ 0\}$. 
\end{theorem}

{\bf The  bi-quadratic algebras on 3 generators.} 
Let $K$  be a field and  $A=K[x_1,x_2,x_3;Q,\mA ,\mB]$ be a bi-quadratic algebra  where $Q=(q_1,q_2,q_3)\in K^{\times 3},$
%\marginpar{AAabc}
\begin{equation}\label{AAabc}
\mA =
 \begin{pmatrix}
  a & b &  c\\
  \alpha & \beta & \g\\
  \lambda & \mu &  \nu\\
  
\end{pmatrix}
\end{equation}
and $\mB=(b_1,b_2,b_3).$  So, the algebra  $A$ is  an algebra that is  generated over the field $K$ by the elements $x_1,x_2$ and $x_3$  subject to  the defining relations:\\
%\marginpar{DRel1}
\begin{equation} \label{DRel1}
x_2x_1-q_1x_1x_2=ax_1+bx_2+cx_3+b_1,
\end{equation} 
%\marginpar{DRel2}
\begin{equation}\label{DRel2}
x_3x_1-q_2x_1x_3=\alpha x_1+\beta x_2+\g x_3+b_2,
\end{equation}
%\marginpar{DRel3}
\begin{equation}\label{DRel3}
x_3x_2-q_3x_2x_3=\lambda x_1+\mu x_2+\nu x_3+b_3.
\end{equation}

{\bf Examples  of bi-quadratic algebras on 3 generators.}\\

 1. The universal enveloping algebra of any  3-dimensional Lie algebra.\\

2. The  {\em 3-dimensional  quantum space} $\mA_{q_1, q_2, q_3}^3:=K [x_1,x_2,x_3;Q,\mA =0 ,\mB =0]$. \\

3. The algebra $U'_q(so_3)$ is generated over the field $K$ by elements $ I_1$, $I_2$ and $ I_3$ subject to the defining relations:
$$q^\frac{1}{2}I_1I_2  - q^{-\frac{1}{2}}I_2I_1= I_3, \;\; 
q^\frac{1}{2}I_2I_3  - q^{-\frac{1}{2}}I_3I_2= I_1, \;\;  
q^\frac{1}{2}I_3I_1  - q^{-\frac{1}{2}}I_1I_3= I_2,$$ 
where $q\in K\backslash 0, \pm 1\}$, \cite{Odesskii-1986,Fairlie-1990}. \\

4. The Askey-Wilson algebras $AW(3)$ introduced by A. Zhedanov, \cite{Zhedanov-AW(3)}. The algebra $AW(3)$ is generated by three elements $K_0$, $K_1$ and $K_2$ subject to the defining relations:
$$[K_0, K_1 ]_w= K_2 , \;\; [K_2, K_0 ]_w= BK_0+C_1K_1+D_1, \;\;[K_1, K_2 ]_w= BK_1+C_0K_0+D_0,$$
where $B$, $C_0$, $C_1$, $D_0$, $D_1\in K$, $[L,M]_w:= wLM-w^{-1}ML$ and $w\in K^\times$. \\

Theorem \ref{11Apr18} gives explicit  necessary and  sufficient conditions for the defining relations of the algebra $A$ to be consistent (i.e., $A\neq0$) and the algebra $A$ has PBW-basis, i.e., $A=\bigoplus_{ \alpha\in \N^3} K x^\alpha $ where $x^\alpha = x_1^{\alpha_1}x_2^{\alpha_2}x_3^{\alpha_3}$ and $\alpha = (\alpha_1,\alpha_2,\alpha_3)$.\\

\begin{theorem}\label{11Apr18}%\marginpar{11Apr18}

Suppose that the algebra $A$ is  generated  over a field $K$ by the elements $x_1,x_2$ and  $x_3$ that satisfy the defining relations (\ref{DRel1}), (\ref{DRel2}) and (\ref{DRel3}). Then the defining relations are consistent and $A=\bigoplus_{ \alpha\in \N^3} K x^\alpha $ where $x^\alpha = x_1^{\alpha_1}x_2^{\alpha_2}x_3^{\alpha_3}$ iff the following conditions hold:\\
%\marginpar{$(X_1X_2)$}
\begin{equation}\label{$(X_1X_2)$}
(1-q_3) \alpha = (1-q_2 )\mu,\ \ \  \ \ \ \ \ \ \ \ \ \ \ \ \ \ \ \ \ \ \ \ \ \ \ \ \ \ \ \ \ \ \ \ \ \ \ \ \ \ \ \ \ \ \ \ \ \ \ \ \ \ \ \ \ \ \ \ \ \ \ \ \ \ \ \ \ \ \ \ \ \ \ \ \ \ 
\end{equation}
%\marginpar{$(X_1X_3)$}
\begin{equation} \label{$(X_1X_3)$}
(1-q_3) a = (1-q_1 )\nu,\ \ \ \  \ \ \ \ \ \ \ \ \ \ \ \ \ \ \ \ \ \ \ \ \ \ \ \ \ \ \ \ \ \ \ \ \ \ \ \ \ \ \ \ \ \ \ \ \ \ \ \ \ \ \ \ \ \ \ \ \ \ \ \ \ \ \ \ \ \ \ \ \ \ \ \ \
\end{equation}
%\marginpar{$(X_2X_3)$}
\begin{equation} \label{$(X_2X_3)$}
 (1-q_2) b = (1-q_1 )\g,\ \ \ \ \ \ \ \ \ \ \ \ \ \ \ \ \ \ \ \ \ \ \ \ \ \ \ \ \ \ \ \ \ \ \ \ \ \ \ \ \ \ \ \ \ \ \ \ \ \ \ \ \ \ \ \ \ \ \ \ \ \ \ \ \ \ \ \ \ \ \ \ \ \ \ \ \
\end{equation}
%\marginpar{$(X_1X_1)$}
\begin{equation} \label{$(X_1X_1)$}
(1-q_1q_{2})\lambda=0,\ \  \ \ \ \ \ \ \ \ \ \ \ \ \ \ \ \ \ \ \ \ \ \ \ \ \ \ \ \ \ \ \ \ \ \ \ \ \ \ \ \ \ \ \ \ \ \ \ \ \ \ \ \ \ \ \ \ \ \ \ \ \ \ \ \ \ \ \ \ \ \ \ \ \ \ \ \ \ \ \ \ \ 
\end{equation}
%\marginpar{$(X_2X_2)$}
\begin{equation} \label{$(X_2X_2)$}
(q_1-q_3)\beta=0,  \ \ \ \ \ \ \ \ \ \ \ \ \ \ \ \ \ \ \ \ \ \ \ \ \ \ \ \ \ \ \ \ \ \ \ \ \ \ \ \ \ \ \ \ \ \ \ \ \ \ \ \ \ \ \ \ \ \ \ \ \ \ \ \ \ \ \ \ \ \ \ \ \ \ \ \ \ \ \ \ \
\end{equation}
%\marginpar{$(X_3X_3)$}
\begin{equation} \label{$(X_3X_3)$}
(1-q_2q_3)c=0,  \ \ \ \ \ \ \ \ \ \ \ \ \ \ \ \ \ \ \ \ \ \ \ \ \ \ \ \ \ \ \ \ \ \ \ \ \ \ \ \ \ \ \ \ \ \ \ \ \ \ \ \ \ \ \ \ \ \ \ \ \ \ \ \ \ \ \ \ \ \ \ \ \ \ \ \ \ \ \ \ 
\end{equation}
%\marginpar{$(X_1)$}
\begin{equation} \label{$(X_1)$}
 ((1-q_3)\alpha-\mu)a+(b+q_1\g)\lambda-\nu\alpha+(q_1q_2-1)b_3=0,  \ \ \ \ \ \ \ \ \ \ \ \ \ \ \ \ \ \ \ \ \ \ \ \ \ \ \ \ \ \ \ \ \ 
\end{equation}
%\marginpar{$(X_2)$}
\begin{equation} \label{$(X_2)$}
 (a- \nu ) \beta +q_1 \g \mu -q_3 \alpha b +(q_1-q_3)b_2=0, \ \ \ \ \ \ \ \ \ \ \ \ \ \ \ \ \ \ \ \ \ \ \ \ \ \ \ \ \ \ \ \ \ \ \ \ \ \ \ \ \ \ \ \ \ \ \ \ \
\end{equation}
%\marginpar{$(X_3)$}
\begin{equation} \label{$(X_3)$}
(a+(q_1-1)\nu)\g +b \nu- (\mu+q_3\alpha)c+(1-q_2q_3)b_1=0,\ \ \ \ \ \ \\ \ \ \ \ \ \ \ \ \ \ \ \ \ \ \ \ \ \ \ \ \ \ \ \ \  \ \ \ \
\end{equation}
%\marginpar{$(1)$}
\begin{equation} \label{$(1)$}
-(\mu +q_3\alpha)b_1+(a-\nu)b_2+(b+q_1\g)b_3=0. \ \ \ \ \\ \ \ \ \ \ \ \ \ \ \ \ \ \ \ \ \ \ \ \ \ \ \ \ \ \ \ \ \ \ \ \ \ \ \ \ \ \ \ \ \ \ \ \ \ \ \ \
 \end{equation}

Furthermore, if $A=\bigoplus_{ \alpha\in \N^3} K x^\alpha $ where $x^\alpha = x_1^{\alpha_1}x_2^{\alpha_2}x_3^{\alpha_3}$ then $A=\bigoplus_{\alpha \in \N^3} K x^\alpha_\sigma$ for all $\sigma \in S_3$ where $x_\sigma^\alpha = x_{\sigma(1)}^{\alpha_1}x_{\sigma(2)}^{\alpha_2}x_{\sigma(3)}^{\alpha_3}$.
\end{theorem}

{\em Remark.} In this paper, when we say `a bi-quadratic algebra on 3 generators' we mean `a bi-quadratic algebra on 3 generators {\em with PBW basis}' (i.e. the ones that satisfy Theorem \ref{11Apr18}).\\

{\bf 4 main classes of bi-quadratic algebras $A$.}  In order to classify  bi-quadratic algebra on 3 generators we have to consider the following 4 classes (in view of the $S_3$-action): 
\begin{enumerate}
\item  $q_1=q_2=q_3=1$, 
\item  $q_1\neq 1$, $q_2=q_3=1$,
\item $q_1\neq 1$, $q_2\neq 1$, $q_3=1$, 
\item   $q_1\neq 1$, $q_2\neq 1$, $q_3\neq 1$.
\end{enumerate}

Explicit descriptions of the algebras in each of the four cases are given in Section \ref{LIETYPE},  Section \ref{N1E1E1}, Section \ref{N1N1E1} and  Section \ref{N1N1N1}, respectively. \\

Till the end of this section,   $$A=K[x_1,x_2,x_3;Q,\mA ,\mB]$$ is  a bi-quadratic algebra on 3-generators where $Q=(q_1,q_2,q_3)\in K^{\times 3}$ and (\ref{AAabc})--(\ref{DRel3}) hold. \\

{\bf 1.  Classification (up to isomorphism) of bi-quadratic algebras on 3 generators of Lie type, i.e. when $q_1=q_2=q_3=1$.} 

\begin{theorem}
\label{9Oct18}%\marginpar{9Oct18}
Let $A$ be the algebra of Lie type. Then the algebra $A$ is isomorphic to one of the following (pairwise non-isomorphic) algebras:
\begin{enumerate}
\item $P_3=K[x_1,x_2,x_3]$, a polynomial algebra in $3$ variables.
\item $U(sl_2(K)),$ the universal enveloping algebra of the Lie algebra $sl_2(K)$.
\item $U(\CH_3),$  the universal enveloping algebra of the Heisenberg  Lie algebra $\CH_3$.
\item $U(\CN)/(c-1)\simeq K\langle x,y,z\ |\ [x,y]=z,\ [x,z]=0, \ [y,z]=1\rangle$  and the algebra $ U(\CN)/(c-1)$  is a tensor product $A_1 \otimes K[x']$ of its subalgebras, the Weyl algebra $A_1=K\langle y,z\ |\ [y,z]=1\rangle$ and the polynomial algebra $K[x']$ where $x'=x+\dfrac{1}{2}z^2$.
\item $U(\gn_2 \times Kz) \simeq K\langle x,y,z \ |\  [x,y]=y\rangle$ and $z$ is a central element.
\item $U(\CM)/(c-1)\simeq K\langle x,y,z\ |\ [x,y]=y,\ [x,z]=1,\ [y,z]=0\rangle$  and the algebra $ U(\CM)/(c-1)$  is a skew polynomial algebra $A_1[y;\sigma] $ where $A_1= K\langle x,z \ |\  [x,z]=1 \rangle$ is the Weyl algebra and $\sigma$ is an automorphism of $A_1$ given by the rule $\sigma(x)=x+1$ and $\sigma(z)=z$.
\end{enumerate} 

\end{theorem}

The proof of Theorem \ref{9Oct18} is given in Section \ref{LIETYPE}. \\

{\bf 2.  Description  of bi-quadratic algebras on 3 generators  when $q_1\neq 1$, $q_2=q_3=1$.}

\begin{theorem}\label{A28Oct18}%\marginpar{A28Oct18}
 Suppose  that $q_1\neq 1$, $q_2=q_3=1$ and  $A=K[x_1,x_2,x_3;Q,\mA ,\mB]$. 
\begin{enumerate}
\item If $\mu +\alpha \neq 0$ then (up to $G_3$):
$$x_2x_1=q_1x_1x_2, \;\; x_3x_1= x_1(x_3+\alpha ) \;\; {\rm and}\;\; x_3x_2= x_2(x_3+\mu )$$
where $(\alpha , \mu ) \in \{ (1,\mu'), (0,1)\, | \, \mu'\in K\backslash \{ -1\} \}$. The elements $x_1$ and $x_2$ are regular normal elements of $A$. The algebra $A=K[x_3][x_2, x_1; \s , \tau , b=0, \rho =q_1]$ is a diskew polynomial algebra where $\tau (x_3) = x_3-\alpha$ and $ \s (x_3) = x_3-\mu$. Furthermore, the algebra $A= \CD [x_2, x_1; \s , \tau , a=h]$ is a GWA where $\CD = K[x_3][h; \tau \s]$ is a skew polynomial ring where $\tau \s (x_3)= x_3-\alpha - \mu$, $\s (h) = q_1 h$ and $\tau (h)= q_1^{-1}h$. The element $h$ is a regular normal element of $A$. 
\begin{enumerate}
\item The algebra $A=K[x_3][x_1;\s_1][x_2;\s_2]$ is an iterated skew polynomial algebra where $\s_1(x_3)=x_3-\alpha$, $\s_2 (x_3)=x_3-\mu$ and $\s_2(x_1)=q_1x_1$. 
\item The algebra homomorphism $A(q_1, \alpha =1, \mu =0)\ra A(q_1^{-1}, \alpha =0, \mu =1)$, $x_1\mapsto x_2$,  $x_2 \mapsto x_1$, $x_3 \mapsto x_3$ is an algebra isomorphism. 
\item In particular, the algebra homomorphism $A(q_1=-1, \alpha =1, \mu =0)\ra A(q_1=-1, \alpha =0, \mu =1)$, $x_1\mapsto x_2$,  $x_2 \mapsto x_1$, $x_3 \mapsto x_3$ is an algebra isomorphism. 
\end{enumerate}

\item If $\mu +\alpha = 0$ then (up to $G_3$):
$$x_2x_1=q_1x_1x_2+cx_3+b_1, \;\; x_3x_1= x_1(x_3+\alpha ) \;\; {\rm and}\;\; x_3x_2= x_2(x_3-\alpha )$$
where exactly one of the following cases occurs: 
\begin{enumerate}
\item $\alpha =0$ and  $c,b_1\in \{ 0,1 \}$, 
\item $\alpha =1$ and   $(c,b_1) \in \{ (0,0),  (1,1)\}$, 
\item  $\alpha =1$, $c=0$, $b_1=1$ or  $\alpha =1$, $c=1$,  $b_1\in K\backslash \{ 1\}$. 
\end{enumerate}
So, the cases (a)--(c) can be written as $(\alpha , c,  b_1)\in \{ 0,1\}^3$ or $\alpha =1$, $c=1$, $b_1\in K\backslash \{0,  1\}$.  The algebra $A=K[x_3][x_2, x_1; \s , \s^{-1} , b=cx_3+b_1, \rho =q_1]$ is a diskew polynomial algebra where $\s (x_3) = x_3+\alpha$. Furthermore, the algebra $A= K[x_3, h] [x_2, x_1; \s , a=h]$ is a classical GWA where $\s (x_3)= x_3+\alpha $ and $\s (h) = q_1h+cx_3+b_1 $.     
\end{enumerate}
\end{theorem}

The proof of Theorem \ref{A28Oct18} is given in Section
 \ref{N1E1E1}. In Section \ref{CLASN=2}, definitions and some results about  generalized Weyl algebra and diskew polynomial rings are given.   \\

{\bf 3.  Description  of bi-quadratic algebras on 3 generators  when  $q_1\neq 1$, $q_2\neq 1$, $q_3=1$.} 

\begin{theorem}\label{28Oct18}%\marginpar{28Oct18}
 Suppose  that $q_1\neq 1$, $q_2\neq 1$,  $q_3=1$ and let $A=K[x_1,x_2,x_3;Q,\mA ,\mB]$. 
\begin{enumerate}

\item If $1-q_1q_2 \neq 0$ then (up to $G_3$):
$$x_2x_1=q_1x_1x_2, \;\; x_3x_1=q_2 x_1x_3 \;\; {\rm and}\;\; x_3x_2= x_2x_3,$$
and so $A\simeq \mA^3_{(q_1,q_2,1)}$.

\item If $1-q_1q_2 = 0$ then (up to $G_3$):
$$x_2x_1=q_1x_1x_2, \;\; x_3x_1=q_1^{-1} x_1x_3 \;\; {\rm and}\;\; x_3x_2= x_2x_3+\l x_1+b_3$$
where $\l , b_3\in \{ 0,1\}$.
 The algebra $A= K[x_1][x_3, x_2; \s , \s^{-1}, \rho =1, b=\l x_1+b_3]$ is a diskew polynomial algebra  where $\s (x_1) = q_1^{-1} x_1$. The algebra $A=K[x_1, h][x_3, x_2; \s , a=h]$ is a classical GWA where $\s (x_1) = q_1^{-1}x_1$ and $\s ( h) = h+b$. If $b_3=0$ then the element $C=h+\alpha'$ where $\alpha' = \l (1-q_1^{-1})^{-1}x_1$ is a central element of the algebra $A$ and $A= K[x_1, C][x_3, x_2; \s , a = C-\alpha']$ is a classical GWA where $\s ( x_1) = q_1^{-1} x_1$ and $\s (C)=C$. 

\end{enumerate}
\end{theorem}

The proof of Theorem \ref{A28Oct18} is given in Section  \ref{N1N1E1}. \\

{\bf 4.  Description  of bi-quadratic algebras on 3 generators  when $q_1\neq 1$, $q_2\neq 1$, $q_3\neq 1$.} 
 The class  of algebras $A$ is a disjoint  union of subclasses that satisfy the condition that elements of the set $\{ q_1-q_3, 1-q_1q_2, 1-q_2q_3\}$   are either zero or not. Potentially, there are $2\times 2 \times 2=8$ subclasses but in fact, there are only 5 since the condition that $q_1-q_3=0$ implies that $1-q_1q_2= 1-q_2q_3$; and the conditions $1-q_1q_2=0$, $1-q_2q_3=0$ imply that $q_1-q_3=0$. 

\begin{itemize}
\item Case 1: $q_1-q_3=0$, $1-q_1q_2=0$. 
\item Case 2: $q_1-q_3=0$, $1-q_1q_2\neq 0$. 
\item Case 3: $q_1-q_3\neq 0$, $1-q_1q_2=0$, $1-q_2q_3\neq 0$. 
\item Case 4: $q_1-q_3\neq 0$, $1-q_1q_2\neq 0$, $1-q_2q_3=0$. 
\item Case 5: $q_1-q_3\neq 0$, $1-q_1q_2\neq 0$, $1-q_2q_3\neq 0$. 
\end{itemize}

{\bf $\bullet$ Case 1: $q_1-q_3=0$ and $1-q_1q_2=0$, i.e., $q:=q_1=q_2^{-1}=q_3\neq 1$ (the quantum bi-quadratic algebras).} \\

{\it Definition.}  A bi-quadratic algebra $A$ is called the  {\bf  quantum bi-quadratic algebra} if $q:=q_1=q_2^{-1}=q_3\neq 1$. The element $q$ is called the {\em quantum parameter} of $A$. \\

In Theorem \ref{17Mar15}, an  expression like `$\mu\in K^\times / K^{\times n}$' means `$\mu$ is unique up to $K^{\times n}$' where $K^{\times n}:=\{ \xi^n\, | \, \xi \in K^\times \}$ and $m\geq 2$.

\begin{theorem}\label{17Mar15}%\marginpar{17Mar15}
Let $A$ be a quantum bi-quadratic algebra over the field $K$ with quantum parameter $q_1$. Then  (up $G_3'$)  
%\marginpar{XBQA1}
\begin{equation} \label{XBQA1}
x_2x_1=qx_1x_2+cx_3+b_1,
\end{equation} 
%\marginpar{XBQA2}
\begin{equation}\label{XBQA2}
x_3x_1=q^{-1}x_1x_3+\beta x_2+b_2,
\end{equation}
%\marginpar{XBQA3}
\begin{equation}\label{XBQA3}
x_3x_2=qx_2x_3+\lambda x_1+b_3,
\end{equation}
where $q\in \{q_1, q_1^{-1}\}$ and (up $G_3'$) the parameters $c,\beta , \l , b_1, b_2, b_3$ belong precisely to one of the four cases below (that correspond to  Cases 1--4 in the proof):\\

Case 1: $(c,\beta , \l )\in (1, K^\times / K^{\times 2}, K^\times / K^{\times 2})$ and $b_1, b_2, b_3\in K$; \\

Case 2: $\l =0$, $(c,\beta ,b_3 )\in (1, K^\times / K^{\times 4}, 1)\coprod (1, K^\times / K^{\times 2},0)$   and $b_1, b_2, \in K$; \\

Case 3: $\beta =\l =0$, $(c,b_2 ,b_3 )\in (1, 1, K^\times / K^{\times 3})\coprod \{ (1,1,0), (1,0,0)\}$   and $b_1 \in K$; \\

Case 4: $c=\beta =\l =0$, $(b_1,b_2 ,b_3 )\in (1, 1, K^\times / K^{\times 2})\coprod \{ (1,1,0), (1,0,0), (0,0,0)\}$. \\

If, in addition,  $\sqrt{K}\subseteq K$ and $ \sqrt[3]{K}\subseteq K$ then  (up $G_3'$) there are the following 11 types of the quantum bi-quadratic algebras:
\begin{enumerate}
\item $c=\beta =\l= 1$ and $b_1, b_2, b_3\in K$. 
\item $c=\beta =1$, $ \l= 0$ and $b_1, b_2\in K$, $b_3\in \{ 0,1\}$.
\item  $c=1$, $\beta =\l  = 0$  and $b_1\in K$, $ (b_2, b_3)\in \{ (1,1), (1,0), (0,0)\}$.
\item  $c  =  \beta  = \l=  0$   and $(b_1, b_2, b_3)\in \{ (1,1,1), (1,1,0), (1,0,0), (0,0,0)\}$.
\end{enumerate}
\end{theorem}

{\it Remark.} If one applies an odd permutation of $S_3$ then the quantum parameter $q$ becomes $q^{-1}$ but for even permutations of $S_3$ it does not change. \\

{\bf $\bullet$ Case 2: $q_1-q_3=0$ and $1-q_1q_2\neq 0$.}

\begin{theorem}\label{QCase-2}%\marginpar{QCase-2}
Suppose that all $q_i\neq 1$,  $q_1-q_3=0$ and $1-q_1q_2\neq 0$. Then (up to $G_3'$):  
$$
x_2x_1=q_1x_1x_2,\;\;  x_3x_1=q_2x_1x_3+\beta x_2+b_2\;\; {\rm and}\;\; 
x_3x_2=q_1x_2x_3
$$
where $\beta , b_2\in \{ 0,1\}$. 
\end{theorem}

{\bf $\bullet$ Case 3: $q_1-q_3\neq 0$, $1-q_1q_2= 0$ and  $1-q_2q_3\neq 0$.} \\

Notice that the third condition (i.e., $1-q_2q_3\neq 0$) follows from the first two.

\begin{theorem}\label{QCase-3}%\marginpar{QCase-3}
Suppose that all $q_i\neq 1$,  $q_1-q_3\neq 0$ and $1-q_1q_2=0$ ($\Rightarrow $  $1-q_2q_3\neq 0$).  Then (up to $G_3'$):  
$$
x_2x_1=q_1x_1x_2,\;\;  x_3x_1=q_1^{-1}x_1x_3\;\; {\rm and}\;\; 
x_3x_2=q_3x_2x_3+\l x_1+b_3
$$
where $\l , b_3\in \{ 0,1\}$. 
\end{theorem}

{\bf $\bullet$ Case 4: $q_1-q_3\neq 0$, $1-q_1q_2\neq 0$ and  $1-q_2q_3=0$.}

\begin{theorem}\label{QCase-4}%\marginpar{QCase-4}
Suppose that all $q_i\neq 1$, $q_1-q_3\neq 0$, $1-q_1q_2\neq 0$ and  $1-q_2q_3=0$.  Then (up to $G_3'$):  
$$
x_2x_1=q_1x_1x_2+cx_3+b_1,\;\;  x_3x_1=q_2x_1x_3\;\; {\rm and}\;\; 
x_3x_2=q_2^{-1}x_2x_3
$$
where $c , b_1\in \{ 0,1\}$. 
\end{theorem}

{\bf $\bullet$ Case 5: $q_1-q_3\neq 0$, $1-q_1q_2\neq 0$ and  $1-q_2q_3\neq 0$.}

\begin{theorem}\label{QCase-5}%\marginpar{QCase-5}
Suppose that all $q_i\neq 1$, $q_1-q_3\neq 0$, $1-q_1q_2\neq 0$ and  $1-q_2q_3\neq 0$.  Then $\mA =0$ and $\mB =0$, i.e., $A\simeq \mA^3_Q$.  
\end{theorem}

Proofs of the above results are given in  Section \ref{N1N1N1}.

%%%%%%%%%%%%%%%%%% SECTION 2 %%%%%%%%%%%%%%%%%%%%%%%%

\section{Classification of bi-quadratic algebras on 2 generators}\label{CLASN=2}%\marginpar{CLASN=2}

It is an easy exercise  to give a classification (up to isomorphism) of bi-quadratic algebras on 2 generators (Theorem \ref{31Oct18}) for $D=K$. 

{\it Definition.} Let $K$ be a field, $q\in K^\times :=K\backslash \{ 0\}$, and $a,b,c,\in K$. Then the algebra
%\marginpar{Abi-2g}
\begin{equation}\label{Abi-2g}
A= K[x_1, x_2; q, a,b,c]:=   K\langle x_1, x_2 \, | \, x_2x_1=qx_1x_2+ax_1+bx_2+c\rangle
\end{equation}
 is a bi-quadratic algebra on 2 generators.

  The algebra $A= K[x_1][x_2; \s , \d ]$ is a skew polynomial algebra where $\s (x_1) = qx_1+b$ and $\d (x_1) = ax_1+c$. So, the algebra $A$ is a Noetherian domain with PBW basis, $A= \bigoplus_{i,j\in  \N} Kx_1^ix_2^j= \bigoplus_{i,j\in  \N} Kx_2^ix_1^j$. and the scalars $a,b,c$ are arbitrary. 

An element of a ring $R$ is called a {\em regular element} if it is neither a left nor right zero-divisor of $R$. The set of all regular elements of the ring $R$ is denoted by $\CC_R$. An element $a$ of the ring $R$ is called a {\em normal element} if $Ra=aR$.   A $\Z$-homomorphism $f: R\ra R$ is called a {\em locally nilpotent map} if $R=\cup_{i\geq 1} \ker (f^i)$. An element $r$ of $R$ is called a {\em locally nilpotent element} if the inner derivation of $R$, $\ad_r:R\ra R$, $x\mapsto rx-xr$ is a locally nilpotent map. 

For a $K$-algebra $A$, we denote by $\CN (A)$ the monoid of all regular normal elements of $A$, $\overline{\CN} (A):=\CN (A) / K^\times$ (the factor monoid modulo its central subgroup $K^\times$) and $\widetilde{\CN}(A)=\CN (A) / Z(A)^{reg}$ (the factor monoid modulo the relation defined by the submonoid  $Z(A)^{reg}:=\CN (A)\cap Z(A)$ of regular central elements of $A$, i.e., $u\sim v$ if $uZ(A)^{reg}\cap vZ(A)^{reg}\neq \emptyset$. 

\begin{theorem}\label{31Oct18}%\marginpar{31Oct18}
Up to isomorphism, there are only five bi-quadratic algebras on 2 generators:
\begin{enumerate}
\item The polynomial algebra $K[x_1,x_2]$,
\item The (first) Weyl algebra $A_1=\langle x_1, x_2\, | \, x_2x_1-x_1x_2=1\rangle $, 
\item The universal enveloping algebra of the Lie algebra $\gn_2=\langle x_1, x_2\, | \, [x_2, x_1]=x_1\rangle$, $U(\gn_2)=K\langle x_1, x_2\, | \, x_2 x_1-x_1x_2=x_1\rangle$, 
\item The quantum plane $\mA^2_2=\langle x_1, x_2\, | \, x_2x_1=qx_1x_2\rangle $ where $q\in K\backslash \{ 0,1\}$, and  
\item The quantum Weyl algebra $A_1(q)=\langle x_1, x_2\, | \, x_2x_1-qx_1x_2=1\rangle $  where $q\in K\backslash \{ 0,1\}$.
\end{enumerate}
\end{theorem}

{\it Proof}.  Let $A=K[x_1,x_2; q, a,b,c ]$ be a bi-quadratic algebra on 2 generators. \\

{\sc Step 1.} {\em The algebra $A$ is isomorphic to one of the algebras in statements 1--5}: \\

(i) {\em Suppose that $q=1$}. 

If $a=b=c=0$ then $A\simeq K[x_1,x_2]$.

If $a=b=0$ and $c\neq 0$ then $A\simeq A_1$.

If $(a, b)\neq (0,0)$  then $A\simeq U(\gn_2)$. If, say $a\neq 0$, then $[a^{-1}x_2, ax_1+bx_2+c]=ax_1+bx_2+c$, and the result follows.

(ii) {\em  Suppose that $q\neq 1$}. By Lemma \ref{a30Oct18}, we can assume that $a=b=0$, i.e., $x_2x_1= qx_1x_2+c$. Up to the change of the variables  $x_1'=\l x_1$, $x_2'= x_2$ (where $\l \in K^\times$) we can assume that $c\in \{ 0,1\}$. If $c=0$ (resp., $c=1$) then $A\simeq \mA^2_q$ (resp., $A\simeq A_1(q)$). \\

{\sc Step 2.} {\em The algebras  in statements 1--5 are not  isomorphic}:  The polynomial algebra $K[x_1, x_2]$ is the only commutative algebra in statements 1--5. So, it remains to show that the algebras in statements 2--5 are not isomorphic. 

If  char $(K)=0$  then the Weyl algebra $A_1$ is the only {\em simple} algebra in statements 2--5. 
If  char $(K)=p>0$  then the Weyl algebra $A_1$ is the only  algebra in statements 2--5 such that all normal elements are central.  So, it remains to show that the algebras in statements 3--5 are not isomorphic.  Let $U=U(\gn_2)$.

(a) $$Z(U)^{reg}= \begin{cases}
K^\times & \text{if  char $K$=0},\\
K[x_1^p, x_2^p]\backslash \{ 0\} & \text{if  char $K$=p}.\\
\end{cases}$$

$$\CN (U) = \begin{cases}
\{ Z(U)^{reg} x_1^i \, | \, i\geq 0\}& \text{if  char $K$=0},\\
\{ Z(U)^{reg} x_1^i \, | \, i=0, 1, \ldots , p-1\} & \text{if  char $K$=p}.\\
\end{cases}$$

 $$\widetilde{\CN}(U) \simeq \begin{cases}
\N & \text{if  char $K$=0},\\
\Z_p:= \Z / p\Z & \text{if  char $K$=p}.\\
\end{cases}$$
 
 (b) $$ Z(\mA^2(q))^{reg}= \begin{cases}
K^\times&   \text{if $q$ is not a root of 1} ,\\
K[x_1^m, x_2^m]\backslash \{ 0\}& \text{if $q$ is a primitive $m$'th a root of 1}.\\
\end{cases}
$$ 
 
  $$\CN (\mA^2(q)) = \begin{cases}
\{  Z(\mA^2(q))^{reg}  x_1^ix_2^j \, | \, i, j\geq 0\} & \text{if  $q$ is not a root of 1},\\
 \{  Z(\mA^2(q))^{reg}  x_1^ix_2^j \, | \, i, j=0, 1, \ldots , m-1\}& \text{if  $q$ is a primitive $m$'th a root of 1}.\\
\end{cases}
 $$ 

$$\widetilde{\CN}(\mA^2(q))  \simeq \begin{cases}
\N \times \N &   \text{if $q$ is not a root of 1} ,\\
\Z_m \times \Z_m& \text{if $q$ is a primitive $m$'th a root of 1}.\\
\end{cases}
$$   
 (c) $$ Z(A_1(q))^{reg}= \begin{cases}
K^\times&   \text{if $q$ is not a root of 1} ,\\
K[x_1^m, x_2^m]\backslash \{ 0\}& \text{if $q$ is a primitive $m$'th a root of 1}.\\
\end{cases}
$$

  $$\CN (A_1(q))) = \begin{cases}
\{  Z(A_1(q))^{reg} h^i \, | \, i\geq 0\} & \text{if  $q$ is not a root of 1},\\
 \{  Z(A_1(q))^{reg}  h^i \, | \, i=0, 1, \ldots , m-1\}& \text{if  $q$ is a primitive $m$'th a root of 1}.\\
\end{cases}
 $$ 

$$\widetilde{\CN}(A_1(q))  \simeq \begin{cases}
\N &   \text{if $q$ is not a root of 1} ,\\
\Z_m & \text{if $q$ is a primitive $m$'th a root of 1}.\\
\end{cases}
$$  

The element $x_1\in U$ is a non-central, regular, normal, locally nilpotent  element of the algebra $U$. There is no such an element in the algebras $\mA^2(q)$ and $A_1(q)$, see (a), (b) and (c). On the other hand, the algebras  $\mA^2(q)$ and $A_1(q)$ are not isomorphic, see (b) and (c) above.  $\Box $

%%%%%%%%%%%%%%%%%% SECTION 3 %%%%%%%%%%%%%%%%%%%%%%%%

\section{ Bi-quadratic algebras on 3 generators}\label{BIQALN=3}%\marginpar{BIQALN=3}

The aim of this section is to give explicit consistency conditions for   3-generated bi-quadratic algebras with PBW basis (Theorem \ref{11Apr18}). Many algebras from the class BQA(3) of bi-quadratic algebras on 3-generators over $K$ are generalized Weyl algebras and diskew polynomial algebras. In this section, we collect results about these two classes of algebras that are used in the proofs of this paper.  \\

{\bf Proof of Theorem \ref{3Jul16}.} The theorem follows at once from the Diamomd Lemma since for the degree-by-lexicographic ordering  the  abmiguities are $\{ x_kx_jx_i\, | \, \leq i<j<k\leq n\}$  and for each of them there are exactly two different ways how to resolve them, see (\ref{kjiL}) and  (\ref{kjiR}).  $\Box $\\

{\bf Proof of Theorem \ref{6Jul16}.} $(\Rightarrow )$ The algebra $A$ is a Lie algebra $(A, [\cdot , \cdot  ])$ where $[a,b]=ab-ba$. By Theorem \ref{3Jul16}, $A=\bigoplus_{\alpha\in \N^n}Dx^\alpha$, and so the vector space $\CG := \bigoplus_{i=0}^nKx_i$ is a Lie subalgebra of $A$ where $x_0=1$. Now, by Theorem \ref{3Jul16}, $A\simeq D\t_K(U(\CG )/(Z-1))$ where $Z=x_0$. 

$(\Leftarrow )$ This implication is obvious. $\Box $\\

{\bf Proof of Theorem \ref{11Apr18}}. In view of Theorem \ref{3Jul16}, we have to equate the coefficients of the only ambiguity $x_3x_2x_1$:
 \begin{eqnarray*}
(x_3x_2)x_1 &=&(q_3x_2x_3+\lambda x_1+\mu x_2+\nu x_3+b_3)x_1\\
& = &q_3x_2(q_2x_1x_3+\alpha x_1+\beta x_2+\g x_3+b_2)\\
&+& \lambda x_1^2+\mu(q_1x_1x_2+ax_1+bx_2+cx_3+b_1)+\nu (q_2x_1x_3+\alpha x_1+\beta x_2+\g x_3+b_2)+b_3x_1\\
&=& q_3q_2(q_1x_1x_2+ax_1+bx_2+cx_3
+b_1)x_3+
q_3\alpha (q_1x_1x_2+ax_1+bx_2+cx_3+b_1)\\
&+&q_3\beta x_2^2+q_3\g  x_2x_3+q_3b_2x_2\\
&=& q_3q_2q_1x_1x_2x_3+(q_1q_3\alpha +q_1\mu )x_1x_2+(q_2q_3a+q_2\nu)x_1x_3+(q_2q_3b+q_3\g)x_2x_3+\lambda x_1^2\\
&+& q_3\beta x_2^2+q_2q_3cx_3^2+(q_3\alpha a+\mu a+b_3+\nu\alpha)x_1+(q_3\alpha b+q_3b_2+\mu b+\nu \beta)x_2\\
&+&(q_3\alpha c+\mu c+\nu\g+q_2q_3b_1)x_3+q_3\alpha b_1+\mu b_1+\nu b_2.
\end{eqnarray*}

\begin{eqnarray*}
x_3(x_2x_1)&=&x_3(q_1x_1x_2+ax_1+bx_2+cx_3+b_1)\\
&=& q_1
(q_2x_1x_3+\alpha x_1+\beta x_2+\g  x_3+b_2)x_2+b(q_3x_2x_3+ \lambda x_1+ \mu  x_2+\nu x_3+b_3)\\
&+& a(q_2x_1x_3+\alpha x_1+\beta x_2+\g x_3+b_2)+cx_3^2+b_1x_3\\
&=& q_1q_2x_1(q_3x_2x_3+\lambda x_1+\mu x_2+\nu x_3+b_3)+\alpha q_1 x_1 x_2+q_1\beta x_2^2\\
&+& q_1\g(q_3x_2x_3+\lambda x_1+\mu x_2+\nu x_3+b_3)+ q_1b_2x_2\\
&=& q_1q_2q_3x_1x_2x_3+(q_1q_2\mu +\alpha q_1)x_1x_2 
+ (q_1q_2\nu +a q_2)x_1x_3 \\
&+& (q_1q_3\g+ q_3b)x_2x_3+ q_1q_2\lambda x_1^2+q_1\beta x_2^2+cx_3^2
+(q_1q_2b_3+q_1\g\lambda +b\lambda +a\alpha)x_1\\
&+& (q_1\g \mu +q_1b_2+b\mu+a\beta)x_2+(q_1\g \nu +b\nu+a\g+b_1)x_3+q_1\g b_3+bb_3+ab_2.
\end{eqnarray*}

%$ \stackrel{\curvearrowleft}{x_3x_2}x_1 =(q_3x_2x_3+\lambda x_1+\mu x_2+\nu x_3+b_3)x_1 = q_3x_2(q_2x_1x_3+\alpha x_1+\beta x_2+\g x_3+b_2)+\lambda x_1^2+\mu(q_1x_1x_2+ax_1+bx_2+cx_3+b_1)+\nu (q_2x_1x_3+\alpha x_1+\beta x_2+\g x_3+b_2)+b_3x_1=q_3q_2(q_1x_1x_2+ax_1+bx_2+cx_3+b_1)x_3+
%q_3\alpha (q_1x_1x_2+ax_1+bx_2+cx_3+b_1)+q_3\beta x_2^2+q_3\g  x_2x_3+q_3b_2x_2
%=q_3q_2q_1x_1x_2x_3+(q_1q_3\alpha +q_1\mu )x_1x_2+(q_2q_3a+q_2\nu)x_1x_3+(q_2q_3b+q_3\g)x_2x_3+\lambda x_1^2+q_3\beta x_2^2+q_2q_3cx_3^2+(q_3\alpha a+\mu a+b_3+\nu\alpha)x_1+(q_3\alpha b+q_3b_2+\mu b+\nu \beta)x_2+(q_3\alpha c+\mu c+\nu\g+q_2q_3b_1)x_3+q_3\alpha b_1+\mu b_1+\nu b_2$.\\

%$x_3\stackrel{\curvearrowleft}{x_2x_1}=x_3(q_1x_1x_2+ax_1+bx_2+cx_3+b_1)=q_1
%(q_2x_1x_3+\alpha x_1+\beta x_2+\g  x_3+b_2)x_2+b(q_3x_2x_3+ \lambda x_1+ \mu  x_2+\nu x_3+b_3)+a(q_2x_1x_3+\alpha x_1+\beta x_2+\g x_3+b_2)+cx_3^2+b_1x_3=q_1q_2x_1(q_3x_2x_3+\lambda x_1+\mu x_2+\nu x_3+b_3)+\alpha q_1 x_1 x_2+q_1\beta x_2^2+q_1\g(q_3x_2x_3+\lambda x_1+\mu x_2+\nu x_3+b_3)+ q_1b_2x_2=q_1q_2q_3x_1x_2x_3+(q_1q_2\mu +\alpha q_1)x_1x_2+(q_1q_2\nu +a q_2)x_1x_3+(q_1q_3\g+ q_3b)x_2x_3+ q_1q_2\lambda x_1^2+q_1\beta x_2^2+cx_3^2
%+(q_1q_2b_3+q_1\g\lambda +b\lambda +a\alpha)x_1+(q_1\g \mu +q_1b_2+b%\mu+a\beta)x_2+(q_1\g \nu +b\nu+a\g+b_1)x_3+q_1\g b_3+bb_3+ab_2.
%$\\

Now, equating  the coefficients of the monomials $x_1x_2, x_1x_3$ and $x_2x_3$ we obtain  the equations (\ref{$(X_1X_2)$}), (\ref{$(X_1X_3)$}) and (\ref{$(X_2X_3)$}) that are multiplied by the nonzero scalars $q_1,\ q_2$ and $q_3$, respectively. So, these equations are equivalent    to the equations (\ref{$(X_1X_2)$}), (\ref{$(X_1X_3)$}) and (\ref{$(X_2X_3)$}),  respectively. Equating  the coefficients of the monomials $x_1^2,\ x_2^2$ and $x_3^2$ we obtain the equations (\ref{$(X_1X_1)$}), (\ref{$(X_2X_2)$}) and (\ref{$(X_3X_3)$}),  respectively. Finally, equating the coefficients of the elements $x_1,\ x_2, \ x_3$ and $1$, we obtain the equations (\ref{$(X_1)$}), (\ref{$(X_2)$}), (\ref{$(X_3)$}) and (\ref{$(1)$}),  respectively. $\Box$\\

{\bf 4 main classes of bi-quadratic algebras $A$.} 
From now on let $$A=K[x_1,x_2,x_3;Q,\mA ,\mB]$$ be a bi-quadratic algebra on 3-generators where $Q=(q_1,q_2,q_3)\in K^{\times 3}$ and (\ref{AAabc}) -- (\ref{DRel3}) hold. 
 In order to classify  bi-quadratic algebra on 3 generators we have to consider the following 4 classes (in view of the $S_3$-action): 
\begin{enumerate}
\item  $q_1=q_2=q_3=1$, 
\item  $q_1\neq 1$, $q_2=q_3=1$,
\item $q_1\neq 1$, $q_2\neq 1$, $q_3=1$, 
\item   $q_1\neq 1$, $q_2\neq 1$, $q_3\neq 1$.
\end{enumerate}

Explicit descriptions of the algebras in each of the four cases are given in Section \ref{LIETYPE},  Section \ref{N1E1E1}, Section \ref{N1N1E1} and  Section \ref{N1N1N1}, respectively. The next lemma is the starting point of analysis of the cases 2--4. 

\begin{lemma}\label{a30Oct18}%\marginpar{a30Oct18}
Suppose that $q_1\neq 1$. Then by making the change of the variables, $$x_1'=x_1-\frac{b}{1-q_1}, \;\; x_2'=x_2-\frac{a}{1-q_1}, \;\; x_3'=x_3,$$ we can assume that $a=b=0$,  and then the equations (\ref{$(X_1X_3)$}) and (\ref{$(X_2X_3)$}) are equivalent to the equalities $\nu =0$ and $\g =0$, respectively. If, in addition, $q_3=1$ then  (\ref{$(X_2X_2)$}) is equivalent to $\beta =0$. 
\end{lemma}

{\it Proof}. By making the above change of variables, the relation (\ref{DRel1}) is equal to the relation $x_2'x_1'-q_1x_1'x_2'=cx_3'+b+\frac{ab}{1-q_1}$. The rest is obvious.   $\Box $\\

Many 3-generated  bi-quadratic algebras are generalized Weyl algebras and diskew polynomial rings. Below, we collected necessary results that are used in  proofs. \\

{\bf (Classical) generalized Weyl algebras $D(\s , a)$ with central element $a$}.

 {\it Definition}, \cite{Bav-GWA-FA-91}-\cite{Bav-Bav-Ideals-II-93}. Let $D$ be a ring, $\s $ be a ring
automorphism of $D$, $a$ is a {\em central} element of $D$. The
{\bf (classical) generalized Weyl algebra} of rank 1 (GWA, for short) $D(\s , a)=D[x,y; \s , a]$ is a
ring generated by the ring $D$  and two elements $x$ and $y$ that
are subject to the defining relations:
%\marginpar{clGWA}
\begin{equation}\label{clGWA}
 xd=\s (d) x\;\;  {\rm and} \;\;  yd=\s^{-1} (d)y\;\; {\rm for \; all} \;\; d\in D,  \;\;
 yx=a \;\; {\rm and} \;\;   xy=\s (a).
\end{equation}
The ring $D$ is called the {\em base ring} of the GWA. The automorphism $\s$ and the element $a$ are called the {\em
defining automorphism} and the {\em defining element} of the GWA, respectively.

 This is an experimental fact that many popular algebras of small
Gelfand-Kirillov dimension are GWAs: the first Weyl algebra $A_1$
and its quantum analogue, the {\em quantum plane}, the {\em
quantum sphere}, $Usl(2)$, $U_qsl(2)$, the {\em Heizenberg}
algebra and its quantum analogues, the $2\times 2$ quantum
matrices, the {\em Witten's} and {\em Woronowicz's} deformations,
Noetherian down-up algebras, etc.

The generalized Weyl algebras were introduced by myself in 1987  when I was an algebra postgradute student at the Taras Shevchenko National University of Kyiv, the Department of Algebra and Mathematical Logic,  and they were the subject of my PhD ``Generalized Weyl algebras and their representations'' submitted at the end of 1990 (defended at the beginning of 1991).\\

{\bf Generalized Weyl algebras with two endomorphisms and a left normal element $a$.} {\it Definition, \cite{GWA-di-skew,QGWA-SCAlg-DPR}}. Let $D$ be a ring, $\s $ and $\tau$ be ring endomorphisms of  $D$, and an element $a\in D$ be  such that
%%\marginpar{abGWA}
\begin{equation}\label{abGWA}
\tau \s (a) = a, \;\; ad= \tau \s (d) a\;\; {\rm and}\;\; \s (a) d= \s\tau (d) \s (a) \;\; {\rm for \; all}\;\; d\in D.
\end{equation}
The {\bf generalized Weyl algebra} (GWA) of rank 1, $A= D(\s, \tau, a) = D[x,y; \s, \tau, a]$,  is a ring generated by $D$, $x$ and $y$ subject to the defining relations:
%\marginpar{GWADEF}
\begin{equation}\label{GWADEF}
xd=\s (d) x\;\;  {\rm and} \;\;  yd=\tau (d)y\;\; {\rm for \; all} \;\; d\in D,
 \;\; yx=a \;\; {\rm and} \;\;   xy=\s (a).
\end{equation}
The ring $D$ is called the {\em base ring} of the GWA $A$. The endomorphisms $\s$, $\tau$ and  the element $a$ are called the {\em
defining endomorphisms} and the  {\em defining element} of the GWA $A$, respectively. By (\ref{abGWA}), the elements $a$ and $\s (a)$ are left normal in $D$. An element $d$ of a ring $D$ is called {\em left normal} (resp., {\em normal}) if $dD\subseteq Dd$ (resp., $Dd=dD$).  To distinguish `old' GWAs from the `new' ones the former are called the {\em classical} GWAs. Every classical GWA is a GWA as the conditions in (\ref{abGWA}) trivially hold if $a$ is central and $\tau = \s^{-1}$.\\

{\bf Diskew polynomial rings}. {\it Definition, \cite{GWA-di-skew,QGWA-SCAlg-DPR}}. Let $D$ be a ring, $\s$ and $\tau$  be its ring endomorphisms, $\rho$ and $b$  be elements of $D$ such that, for all $d\in D$,
%\marginpar{wbd2}
\begin{equation}\label{wbd2}
\s\tau (d)\rho = \rho \tau \s (d)\;\; {\rm and}\;\; \s \tau (d)b=bd,
\end{equation}
The {\bf diskew polynomial ring} (DPR)  $E:= D(\s , \tau, b , \rho ):= D[x,y; \s , \tau , b,\rho ]$ is a ring generated by $D$, $x$ and $y$ subject to the defining relations:
%\marginpar{DiskewDEF}
\begin{equation}\label{DiskewDEF}
xd=\s (d) x\;\; {\rm and}\;\; yd= \tau (d) y\;\; {\rm for\; all}\;\; d\in D, \;\; xy-\rho yx = b.
\end{equation}
By (\ref{wbd2}), $b$ is a left normal element of $D$. If $\tau\s$ (resp., $\s\tau$) is an epimorphism then $\rho$ is a left (resp., tight) normal element of $D$. \\

{\bf Diskew polynomial rings are GWAs when $\rho$ is a unit}. If the element $\rho$ is a unit in $D$ then every diskew polynomial ring is a generalized Weyl algebra, Theorem \ref{2Apr}, where $ \o_\rho : D\ra D$, $d\mapsto \rho d \rho^{-1}$ is the inner automorphism of $D$  that is determined by the unit  $\rho$.
\begin{theorem}\label{2Apr}%\marginpar{2Apr}
(\cite{GWA-di-skew,QGWA-SCAlg-DPR}) Let $E=D[x,y; \s , \tau , b, \rho ]$ be a diskew polynomial ring. Suppose that $\rho$ is a unit in $D$. Then   $x$ and $y$  are  left regular elements of $E$ and  the ring $E= \CD [ x,y; \s , \tau , a=h]$ is a GWA with base ring $\CD := D[h;  \tau \s]$ which is a  skew polynomial ring, $\s$ and $\tau$ are ring endomorphisms of $\CD$ that are extensions of the ring endomorphisms $\s$ and $\tau$ of $D$, respectively, defined  by the rule $\s (h) = \rho h +b$ and $\tau (h) = \tau (\rho^{-1}) (h-\tau (b))$. In particular, $\tau \s (h)=h$ and $\s\tau (h) = \o_\rho (h) = \rho \tau \s (\rho^{-1})h$. Furthermore, $\s\tau = \o_\rho \tau \s$ in $\CD$.
\end{theorem}
%Theorem \ref{2Apr} is a generalization of a similar result for rings $D\langle \s ; b, \rho \rangle $, {\cite[Lemma 1.2, Corollary 1.4]{Bav-GlGWA-1996}}.

{\bf The canonical left normal element $C$ of a diskew polynomial ring}.
 Theorem \ref{B6Apr} is a criterion for an element $C=h+\alpha$ (where $\alpha \in D$) to be a left normal element in $E$. Theorem \ref{B6Apr} is a generalization of {\cite[Lemma 1.3]{Bav-GlGWA-1996}}.

\begin{theorem}\label{B6Apr}%\marginpar{B6Apr}
(\cite{GWA-di-skew,QGWA-SCAlg-DPR}) Let $E= D[x,y; \s , \tau , b, \rho ]$ be a diskew polynomial ring such that $\rho$ is a unit, $\CD = D[h; \nu = \tau \s ]$ and $C= h+\alpha$ where $h=yx$ and $\alpha \in D$.  The following statements are equivalent.
\begin{enumerate}
\item The element $C$ is left normal in $E$.
\item $\rho \alpha - \s (\alpha ) = b$, $\nu (\alpha )=\alpha $  and $\alpha d= \nu (d) \alpha$ for all elements $d\in D$.

If one of the equivalent conditions holds then $[h,C]=0$ and
\begin{enumerate}
\item  $C=\rho^{-1} (xy+\s (\alpha ))$, $xC= \rho Cx$ and $yC= \tau (\rho^{-1}) C y$.
\item $E\simeq D[C; \nu ] [x,y; \s , \tau , a:= C-\alpha ]$ is a GWA where $\s (C) = \rho  C$ and $\tau (C) = \tau (\rho^{-1}) C$.
\item The element $C$ is a left normal, left regular element of $E$ and
 $E/ (C)\simeq D [x,y; \s , \tau , -\alpha ]$ is a GWA.
\item The element $C$ is a normal element in $E$ iff $\im (\nu ) = D$.
\item The element $C$ is regular  iff  $C$ is  right regular iff $\ker (\nu  ) =0$.
\item The element $C$ is a normal,  regular element iff $\nu$ is an automorphism of $D$.
\end{enumerate}
\end{enumerate}
\end{theorem}

A more general situation when $E$ is a GWA is described below. 

\begin{theorem}\label{A6Apr}%\marginpar{A6Apr}
(\cite{GWA-di-skew,QGWA-SCAlg-DPR}) Let $E= D[x,y; \s , \tau , b, \rho ]$ be a diskew polynomial ring such that $\rho$ is a unit and $\CD = D[h; \nu = \tau \s ]$ where $h=yx$.  The following statements are equivalent.
\begin{enumerate}
\item There exists an element $C= h+\alpha \in \CD $, where $\alpha \in D$, such that $Cd=\nu (d) C$ for all elements $d\in D$ and $\s (C) = \rho C$.
\item There is an element $\alpha \in D$ such that $\rho \alpha - \s (\alpha ) = b$  and $\alpha d= \nu (d) \alpha$ for all elements $d\in D$.

If one of the equivalent conditions holds then $[h,C]=(\nu (\alpha ) - \alpha ) C$ and
\begin{enumerate}
\item  $C=\rho^{-1} (xy+\s (\alpha ))$, $xC= \rho Cx$ and $yC= \tau (\rho^{-1}) (C+\nu (\alpha ) -\alpha ) y$.
\item $E\simeq D[C; \nu ] [x,y; \s , \tau , a:= C-\alpha ]$ is a GWA where $\s (C) = \rho  C$ and $\tau (C) = \tau (\rho^{-1}) (C+\nu (\alpha ) -\alpha ) $. Furthermore, $\tau \s (C) = C+\nu (\alpha ) - \alpha$ and $\s \tau (C)=\s \tau (\rho^{-1})(\rho C+\s (\nu (\alpha ) - \alpha ))$.
\end{enumerate}
\end{enumerate}
\end{theorem}

{\bf The canonical central  element $C$ of a diskew polynomial ring (under certain conditions)}.
The next corollary is a criterion for an element $C=h+\alpha$ (where $\alpha \in D$) to be a central element in $E$. It follows straightaway from Theorem \ref{B6Apr}. This is a generalization of a similar result for the rings $D\langle \s, b , \rho \rangle$,  and {\cite[Lemma 1.5]{Bav-GlGWA-1996}}.

\begin{corollary}\label{aB6Apr}%\marginpar{aB6Apr}
(\cite{GWA-di-skew,QGWA-SCAlg-DPR})  Let $E= D[x,y; \s , \tau , b, \rho ]$ be a diskew polynomial ring such that $\rho$ is a unit, $\CD = D[h; \nu = \tau \s ]$ and $C= h+\alpha$ where $h=yx$ and $\alpha \in D$.  The following statements are equivalent.
\begin{enumerate}
\item The element $C$ is a central element of $E$.

\item $\rho =1$, $\nu =1$, $\alpha - \s (\alpha ) = b$, and the element $\alpha$  belongs to the centre of $D$.

If one of the equivalent conditions holds then
\begin{enumerate}
\item  $C=xy+\s (\alpha )$.
\item $E\simeq D[C] [x,y; \s , \tau , a:= C-\alpha ]$ is a GWA where $\s (C) = C$ and $\tau (C)=C$.
\item The element $C$ is a  regular element of $E$.
\end{enumerate}
\end{enumerate}
\end{corollary}

%%%%%%%%%%%%%%%%%%   Section  4  %%%%%%%%%%%%%%%%%%%%%%%%%%

\section{The bi-quadratic algebras with $q_1\neq 1$, $q_2=q_3=1$}\label{N1E1E1}%\marginpar{N1E1E1}

The aim of this section is to prove Theorem \ref{A28Oct18} that gives an explicit description of  bi-quadratic algebras with $q_1\neq 1$, $q_2=q_3=1$. \\

{\bf Proof of Theorem \ref{A28Oct18} }.  Since $q_1\neq 1$ and $q_3=1$,   we can assume that $a=b=\nu = \g = \beta =0$, by Lemma \ref{a30Oct18}.  Since $q_1\neq 1$ and $q_2=1$, we must have $1-q_1q_2=1-q_1\neq 0$. Then the conditions (\ref{$(X_1X_2)$}),  (\ref{$(X_1X_1)$}) and (\ref{$(X_3X_3)$})  can be written, respectively,  as follows
$$ 0=0, \;\;\l =0 \;\; {\rm and }\;\; 0=0.$$
Then the conditions (\ref{$(X_1)$})--(\ref{$(1)$})  can be written, respectively,  as follows: $b_3=0$, $ b_2=0$, 
%\marginpar{macb}
\begin{equation}\label{macb}
(\mu +\alpha ) c=0 \;\; {\rm and}\;\; (\mu +\alpha )b_1=0.
\end{equation}
Then the defining relations of the algebra $A$ take the form
\begin{eqnarray*}
x_2x_1 &=& q_1x_1x_2+cx_3+b_1,  \\
x_3x_1 &=& x_1(x_3+\alpha ), \\
x_3x_2 &=& x_2(x_3+\mu ),
\end{eqnarray*}
 where the elements $c,b_1, \alpha$ and $\mu$ satisfy (\ref{macb}).\\
 
1. {\it Suppose that $\mu +\alpha \neq 0$}.  Then the conditions in (\ref{macb}) are equivalent to $c=b_1=0$, and the defining relations of $A$ take the form 
$$x_2x_1=q_1x_1x_2, \;\; x_3x_1= x_1(x_3+\alpha ) \;\; {\rm and}\;\; x_3x_2= x_2(x_3+\mu )$$
where $\alpha , \mu \in K$. By making the change of variables $x_1'=x_1$, $x_2'=x_2$, $x_3'=\o x_3$ (where $\o \in K^\times $), we can assume that 
 $(\alpha , \mu ) \in \{ (1,\mu'), (0,1)\, | \, \mu'\in K\backslash \{ -1\} \}$ and the first part of statement 1 follows.  The second part follows from Theorem \ref{2Apr}.

 The statements (a)--(c) are  obvious.  

2. {\it Suppose that $\mu +\alpha =0$}. Then $\mu = -\alpha$, the conditions in (\ref{macb}) hold automatically, and the defining relations of the algebra $A$ are as follows 
$$x_2x_1=q_1x_1x_2+cx_3+b_1, \;\; x_3x_1= x_1(x_3+\alpha ) \;\; {\rm and}\;\; x_3x_2= x_2(x_3-\alpha )$$
where $c, \alpha , b_1\in  K$. There are three case: 

(a) $\alpha =0$, 

(b) $\alpha \neq 0$ and $ c\alpha = b_1$, and

(c)  $\alpha \neq 0$ and $ c\alpha \neq b_1$. 

The cases (a)--(c) correspond to the cases (a)--(c) of statement 2, respectively. 

(a) If $\alpha =0$, $c=0$ and $b_1\neq 0$ then by making the change of the variables $x_1'=b_1^{-1}x_1$, $x_2'= x_2$, $x_3'=x_3$ the triple $(\alpha, c, b_1)=(0,0,b_1)$ is transformed into the triple $(0,0,1)$.  

If $\alpha =0$, $c\neq 0$ and $b_1\neq 0$ (resp., $b_1=0$ ) then by making the change of the variables $x_1'=b_1^{-1}x_1$, $x_2'= x_2$, $x_3'=b_1^{-1}cx_3$ (resp.,  $x_1'=x_1$, $x_2'= x_2$, $x_3'=cx_3$ ) the triple $(\alpha, c, b_1)=(0,c,b_1)$ is transformed into the triple $(0,1,1)$ (resp., $(0,1,0)$).

(b) Suppose that $\alpha \neq 0$ and $ c\alpha = b_1$. There are two cases:

(b1) $c=0$ and $b_1=0$, and

(b2) $c\neq 0$ and $b_1\neq 0$.

In the case (b1) (resp., (b2)), the substitution $x_1'=x_1$, $x_2'= x_2$, $x_3'=\alpha^{-1}x_3$ (resp., $x_1'=(c\alpha )^{-1}x_1$, $x_2'= x_2$, $x_3'=\alpha^{-1}x_3$) transforms the triple $(\alpha, c, b_1)$ into the triple $(1,0,0)$ (resp., $(1,1,1)$).

(c) Suppose that $\alpha \neq 0$ and $ c\alpha \neq  b_1$. There are two cases:

(c1) $c=0$, and

(c2) $c\neq 0$.

If $c=0$ then $b_1\neq 0$ and the substitution $x_1'=b_1^{-1}x_1$, $x_2'= x_2$, $x_3'=\alpha^{-1}x_3$  transforms the triple $(\alpha, c=0, b_1)$ into the triple $(1,0,1)$.

If $c\neq0$ then $b_1\neq c\alpha$ and the substitution $x_1'=(c\alpha)^{-1}x_1$, $x_2'= x_2$, $x_3'=\alpha^{-1}x_3$  transforms the triple $(\alpha, c, b_1)$ into the triple $(1,1, b_1)$ where $b_1\in K\backslash \{ 1\}$.

So, the first part of statement 2 follows. The second part follows from Theorem \ref{2Apr}. $\Box $\\

%**********  $\gldim +  \Kdim  $   *******

%%%%%%%%%%%%%%%%%%   Section  5  %%%%%%%%%%%%%%%%%%%%%%%%%%

\section{The bi-quadratic algebras with $q_1\neq 1$, $q_2\neq 1$, $q_3=1$}\label{N1N1E1}%\marginpar{N1N1E1}

The aim of this section is to prove Theorem \ref{28Oct18} that gives an explicit description of   bi-quadratic algebras with $q_1\neq 1$, $q_2\neq 1$ and $q_3=1$. \\

{\bf Proof of Theorem \ref{28Oct18}.} Since $q_1\neq 1$ and $q_3=1$,  we can assume that $a=b=\nu = \g = \beta =0$, by Lemma \ref{a30Oct18}. Then the equalities (\ref{$(X_1X_2)$}),
 (\ref{$(X_3X_3)$}), (\ref{$(X_2)$}), (\ref{$(X_3)$}) and (\ref{$(1)$}) can be written respectively as 
$$ \mu =0, \;\; c=0,\;\; b_2=0, \;\; b_1=0\;\; {\rm and}\;\; 0=0.$$
Then the equalities  (\ref{$(X_1X_1)$}) and  (\ref{$(X_1)$})
 are equal respectively to the equalities
$$ (1-q_1q_2)\l =0\;\; {\rm and}\;\; (1-q_1q_2)b_3 =0,$$
and the defining relations are 
$$ x_2x_1= q_1x_1x_2, \;\; x_3x_1= x_1(q_2x_3+\alpha ) \;\; {\rm and } \;\; x_3x_2= x_2x_3+\l x_1 +b_3.$$
Using $x_3'=x_3-\frac{\alpha}{1-q_2}$, we can assume that $\alpha =0$. \\

1. {\em Suppose that} $1-q_1q_2\neq 0$. Then $\l= b_3=0$ and statement 1 follows.\\

2. {\em Suppose that} $1-q_1q_2=0$. Then $\l$ and $b_3$ are arbitrary scalars. Up to change of the variables  $x_1'=\l_1x_1$, $x_2'=\l_2x_2$, $x_3'=x_3$  (where $\l_i\in K^\times$), the elements $\l$ and $b_3$ belong to the set $\{ 0,1\}$ and the first part of statement 2 follows.\\

Clearly, the algebra $A$ is the diskew polynomial algebra as in statement 2. If $b_3=0$ then the element $\alpha' = \l (1-q_1^{-1})^{-1}x_1$ is a solution to the equation $\alpha' - \s (\alpha') = b$ (where $b = \l x_1$). By  Corollary \ref{aB6Apr}, the element $C= h-\alpha'$ is a central element of $A$ and the algebra $A$ is the GWA as in statement 2. 
 $\Box$

%***********  $\gldim +  \Kdim  $   *******

%%%%%%%%%%%%%%%%%%   Section  6  %%%%%%%%%%%%%%%%%%%%%%%%%%

\section{The bi-quadratic algebras with $q_1\neq 1$, $q_2\neq 1$, $q_3\neq 1$}\label{N1N1N1}%\marginpar{N1N1N1}

 In this section we assume that $q_1\neq 1$, $q_2\neq 1$, $q_3\neq 1$.
So, this class of algebras is closed under the action of the symmetric group $S_3$ (under permutation of the canonical generators $x_1, x_2, x_3$). We split the class  of algebras into a disjoint union of subclasses based on the condition that elements of the set $\{ q_1-q_3, 1-q_1q_2, 1-q_2q_3\}$   are  either zero or not. Potentially, there are $2\times 2 \times 2=8$ subclasses but in fact, there are only 5 since the condition that $q_1-q_3=0$ implies  that $1-q_1q_2= 1-q_2q_3$; and the conditions $1-q_1q_2=0$ and  $1-q_2q_3=0$ imply that $q_1-q_3=0$. 

\begin{itemize}
\item Case 1: $q_1-q_3=0$, $1-q_1q_2=0$. 
\item Case 2: $q_1-q_3=0$, $1-q_1q_2\neq 0$. 
\item Case 3: $q_1-q_3\neq 0$, $1-q_1q_2=0$, $1-q_2q_3\neq 0$. 
\item Case 4: $q_1-q_3\neq 0$, $1-q_1q_2\neq 0$, $1-q_2q_3=0$. 
\item Case 5: $q_1-q_3\neq 0$, $1-q_1q_2\neq 0$, $1-q_2q_3\neq 0$. 
\end{itemize}

The next lemma is often used in proofs. 

\begin{lemma}
\label{b15Oct18}%\marginpar{b15Oct18}
Suppose that $R$ is a $K$-algebra,  $q_1,q_2 \in K$ and $q_1 \neq1$.   When we make the substitution $x=x'+\dfrac{1}{1-q_1}$ into the equalities below that hold in the algebra $R$ (where $x,y,z,z' , u,u'\in R$): 

$yx=(q_1x+1)y+x+\dfrac{1}{q_1-1},\  zx=\Bigg(q_2x+\dfrac{1-q_2}{1-q_1}\Bigg)z,\  xz'=z'\Bigg(q_2x+\dfrac{1-q_2}{1-q_1}\Bigg), \ ux=(q_1x+1)u$, and $xu'=u'(q_1x+1),$ we obtain  respectively the equalities:  $$yx'=x'(q_1y+1),\ zx'=q_2x'z,\ x'z'=q_2z'x',\ ux'=q_1x'u\;\; {\rm  and} \;\;  x'u'=q_1u'x'.$$
\end{lemma}

{\it Proof.} The first equality follows from $$yx'+\dfrac{y}{1-q_1}=yx=(q_1x+1)y+x+\dfrac{1}{q_1-1} = \Bigg(q_1x'+\dfrac{q_1}{1-q_1}+1\Bigg)y+x_1'=x_1'(q_1y+1)+\dfrac{y}{1-q_1}.$$ The second equality follows from $$zx'+\dfrac{z}{1-q_1}=zx=\Bigg(q_2x+\dfrac{1-q_2}{1-q_1}\Bigg)z=\Bigg(q_2x'+\dfrac{q_2}{1-q_1}+\dfrac{1-q_2}{1-q_1}\Bigg)z= q_2x'z+\dfrac{z}{1-q_1}.$$ The third equality can be proven in the same manner  as the second.  The fourth equality follows from $$ux'+\dfrac{u}{1-q_1}=ux=(q_1x+1)u = \Bigg(q_1x'+\dfrac{q_1}{1-q_1}+1\Bigg)u=q_1x'u+\dfrac{u}{1-q_1}.$$  The fifth equality is proven in a similar way as the fourth. $\Box$\\

Since $q_1\neq 1$, we can assume that $a=b=\nu =\g =0$, by Lemma \ref{a30Oct18}. Then the equation (\ref{$(X_1X_2)$}) can be written as 
$$ \mu = \frac{1-q_3}{1-q_2}\alpha .$$
In particular, the defining relations of the algebra $A$ are 
$$x_2x_1=q_1x_1x_2+cx_3+b_1,\;\; 
x_3x_1=x_1(q_2x_3+\alpha )+\beta  x_2+b_2, \;\; 
x_3x_2=x_2\Bigg( q_3x_3+\frac{1-q_3}{1-q_2}\alpha\Bigg) +\l x_1+b_3.$$
If $\alpha \neq 0$ then using Lemma \ref{b15Oct18} we can assume that $\alpha =0$. In more detail, first by making the substitution $x_3'=\alpha^{-1} x_3$, we can assume that $\alpha =1$. Then applying 
Lemma \ref{b15Oct18} where $x_3=x_3'+\frac{1}{1-q_2}$ we obtain the result. \\

Summarizing, we can assume that 
%\marginpar{abna}
\begin{equation}\label{abna}
a=b=\alpha = \g = \mu = \nu =0,
\end{equation}
the equalities (\ref{$(X_1X_2)$})--(\ref{$(X_2X_3)$}) automatically hold, and the defining relations are 
%\marginpar{abna1}
\begin{equation} \label{abna1}
x_2x_1=q_1x_1x_2+cx_3+b_1,
\end{equation} 
%\marginpar{abna2}
\begin{equation}\label{abna2}
x_3x_1=q_2x_1x_3+\beta x_2 +b_2,
\end{equation}
%\marginpar{abna3}
\begin{equation}\label{abna3}
x_3x_2=q_3x_2x_3+\lambda x_1+b_3.
\end{equation}

\subsection{Case 1: $q_1-q_3=0$ and $1-q_1q_2=0$, i.e., $q:=q_1=q_2^{-1}=q_3\neq 1$ (the quantum bi-quadratic algebras)}

{\it Definition.}  A bi-quadratic algebra $A$ is called the  {\bf  quantum bi-quadratic algebra} if $q:=q_1=q_2^{-1}=q_3\neq 1$. The element $q$ is called the {\em quantum parameter} of $A$. \\

In particular all elements of the set $\{ q_1-q_3, 1-q_1q_2=1-q_2q_3\}$   are  equal to $0$.  Then {\em the equalities (\ref{$(X_1X_1)$}) --(\ref{$(1)$}) hold automatically}. So, (\ref{abna1})--(\ref{abna3}) are the defining relations for the quantum bi-quadratic algebra $A$ where there is no restrictions on the defining constants  apart from $q_1=q_2^{-1}=q_3\neq 1$. So, the defining relations can be written as in Theorem \ref{17Mar15}. \\

{\bf Classification of orbits of four actions of $\mT^3$ on $K^3$.}
 The 3-dimensional algebraic torus $\mT^3=K^{\times 3}$ acts in four different ways on the affine space $K^3$: For all $\l = (\l_1, \l_2, \l_3)\in \mT^3$ and $\xi = (\xi_1, \xi_2, \xi_3)\in K^3$, \\
 
 {\em Case 1:} $\l \cdot \xi= (\frac{\l_3}{\l_1 \l_2}\xi_1, \frac{\l_2}{\l_1 \l_3}\xi_2, \frac{\l_1}{\l_2 \l_3}\xi_3)$. \\

{\em Case 2:} $\l \cdot \xi= (\frac{\l_3}{\l_1 \l_2}\xi_1, \frac{\l_2}{\l_1 \l_3}\xi_2, \frac{1}{\l_2 \l_3}\xi_3)$. \\

{\em Case 3:} $\l \cdot \xi= (\frac{\l_3}{\l_1 \l_2}\xi_1, \frac{1}{\l_1 \l_3}\xi_2, \frac{1}{\l_2 \l_3}\xi_3)$.\\

{\em Case 4:} $\l \cdot \xi= (\frac{1}{\l_1 \l_2}\xi_1, \frac{1}{\l_1 \l_3}\xi_2, \frac{1}{\l_2 \l_3}\xi_3)$.\\

 The four actions and the classifications of theirs orbits (Proposition \ref{XA17Mar15})  are an essential part of the proof of  Theorem \ref{17Mar15}. The set
 $$K^3= K^{\times 3}\coprod K^3_{\rm sing}$$ is a disjoint union of two $\mT^3$-invariant  sets where $K^3_{\rm sing}:= K^3\backslash  K^{\times 3}= \{ (\xi_1, \xi_2, \xi_3 )\in K^3 \, | \, \xi_1 \xi_2 \xi_3=0\}$. For each element $\xi = (\xi_1, \xi_2, \xi_3 )\in K^3$, the elements $\supp (\xi ) :=(\varepsilon_1, \varepsilon_2, \varepsilon_3)\in \{ 0,1\}^3$ is called the {\em support} of $\xi$ where 
 $$\varepsilon_i=\begin{cases}
1& \text{if } \xi_i\neq 0,\\
0& \text{if } \xi_i= 0.\\
\end{cases}$$

Notice that 
$$ K^3_{\rm sing} =\coprod_{i=1}^3 K^3_{{\rm sing}, i}$$
is a disjoint union of $\mT^3$-invariant subsets $K^3_{{\rm sing}, i}:=\{ \xi \in K^3_{\rm sing}\, | \, z(\xi ) = i\}$ where $z(\xi )$ is the number of zero coordinates of the vector $\xi$. The set 
$$K^3_{{\rm sing}, 1}=\coprod_{i=1}^3K^3_{{\rm sing}, 1}(i)$$
is a disjoint union of $\mT^3$-invariant subsets $K^3_{{\rm sing}, 1}(i):=\{ \xi \in K^3_{{\rm sing},1}\, | \, \xi_i = 0\}$.

For each natural number $n\geq 1$, the image of the group homomorphism $\mu_n : K^\times \ra K^\times$, $a\mapsto a^n$ is denoted by $K^{\times n}$, and $K^\times / K^{\times n}$ is the factor group of the group $K^\times$ modulo $K^{\times n}$. Notice that   $K^\times = K^{\times n}$ for $n=2$ (resp., $n=3$) iff the field $K$ contains all quadratic (resp., cubic) roots of all elements of $K$. 

We denote by $K^3/\mT^3 = \{ \mT^3\xi \, | \, \in K^3\}$ the set of $T^3$-orbits in $K^3$. Clearly, 
  $$K^3/\mT^3 = K^{\times 3}/\mT^3 \coprod K^3_{\rm sing}/\mT^3 =K^{\times 3}/\mT^3 \coprod \coprod_{i=1}^3K^3_{{\rm sing}, i}.$$
  
 In view of the equality above, Proposition \ref{XA17Mar15} classifies orbits of the four actions of $\mT^3$ on $K^3$ (Cases 1--4).

\begin{proposition}\label{XA17Mar15}%\marginpar{XA17Mar15}

\begin{enumerate}
\item For the action $\l \cdot \xi= (\frac{\l_3}{\l_1 \l_2}\xi_1, \frac{\l_2}{\l_1 \l_3}\xi_2, \frac{\l_1}{\l_2 \l_3}\xi_3)$, 
\begin{enumerate}
\item the map
$$ K^{\times 3}/\mT^3 \ra K^\times / K^{\times 2}\times  K^\times / K^{\times 2}, \;\; (\xi_1, \xi_2, \xi_3)\mapsto (\xi_1 \xi_2 K^{\times 2}, \xi_1\xi_3K^{\times 2})$$
is a bijection with the inverse $(\rho K^{\times 2}, \eta K^{\times 2}) \mapsto  \mT^3(1,\rho , \eta )$. 
\item The set $(1, K^\times / K^{\times 2}, 0)\coprod (1, 0, K^\times / K^{\times 2})\coprod (0, 1,  K^\times / K^{\times 2})\coprod \{   (1,0,0), (0,1,0),(0,0,1), (0,0,0)\}$ is a set of representatives of the $\mT^3$-orbits in $K^3_{\rm sing}/\mT^3$. 
\end{enumerate}

\item For the action $\l \cdot \xi= (\frac{\l_3}{\l_1 \l_2}\xi_1, \frac{\l_2}{\l_1 \l_3}\xi_2, \frac{1}{\l_2 \l_3}\xi_3)$,
\begin{enumerate}
\item  the map
$$ K^{\times 3}/\mT^3 \ra K^\times / K^{\times 4}, \;\; (\xi_1, \xi_2, \xi_3)\mapsto \frac{\xi_2}{\xi_1\xi_3^2}K^{\times 4}$$
is a bijection with the inverse $\rho K^{\times 4}\mapsto  \mT^3 (1, \rho , 1)$. 

\item The set $(1, K^\times / K^{\times 2}, 0)\coprod \{   (1,0,1), (0,1,1), (1,0,0), (0,1,0), (0,0,1), (0,0,0)\}$ is a set of representatives of the $\mT^3$-orbits in $K^3_{\rm sing}/\mT^3$. 
\end{enumerate}

\item For the action $\l \cdot \xi= (\frac{\l_3}{\l_1 \l_2}\xi_1, \frac{1}{\l_1 \l_3}\xi_2, \frac{1}{\l_2 \l_3}\xi_3)$,

\begin{enumerate}
\item  the map
$$ K^{\times 3}/\mT^3 \ra K^\times / K^{\times 3}, \;\; (\xi_1, \xi_2, \xi_3)\mapsto \frac{\xi_3}{\xi_1\xi_2^2}K^{\times 3}$$
is a bijection with the inverse $\rho K^{\times 3}\mapsto  \mT^3(1,  1, \rho )$.
\item The set $\{ 0,1\}^3\backslash \{ (1,1,1)\}$  is a set of representatives of the $\mT^3$-orbits in $K^3_{\rm sing}/\mT^3$. 
\end{enumerate}

\item For the action $\l \cdot \xi= (\frac{1}{\l_1 \l_2}\xi_1, \frac{1}{\l_1 \l_3}\xi_2, \frac{1}{\l_2 \l_3}\xi_3)$,
\begin{enumerate}
\item  the map
$$ K^{\times 3}/\mT^3 \ra K^\times / K^{\times 2}, \;\; (\xi_1, \xi_2, \xi_3)\mapsto \frac{\xi_3}{\xi_2\xi_2}K^{\times 2}$$
is a bijection with the inverse $\rho K^{\times 2}\mapsto  \mT^3(1,1, \rho )$. 
\item The set $\{ 0,1\}^3\backslash \{ (1,1,1)\}$  is a set of representatives of the $\mT^3$-orbits in $K^3_{\rm sing}/\mT^3$. 
\end{enumerate}

\item For each of the four actions above,
 the map 
$$( K^{\times 3}_{{\rm sing}, 2} \coprod K^{\times 3}_{{\rm sing}, 3})/\mT^3 \ra\{ (0,0,0),  (1,0,0), (0,1,0),(0,0,1)\}, \;\; \xi\mapsto \supp (\xi)$$
is a bijection with the inverse $v\mapsto \mT^3v$.

\end{enumerate}
\end{proposition}

{\it Proof}.  Given $\xi = (\xi_1, \xi_2, \xi_3)\in K^{\times 3}$ and $\l =(\l_1, \l_2, \l_3)\in \mT^3$.

1(a) $\l \cdot \xi \in (1, K^\times , K^\times )$ iff $\frac{\l_3}{\l_1\l_2}\xi_1=1$, and in this case  $\l \cdot \xi =(1, \xi_1\xi_2\l_1^{-2}, \xi_1\xi_3\l_2^{-2})$; and the statement (a) follows. 

(b) If $\xi\in K^3_{\rm sing}$ then at least one of the coordinates of the vector $\xi$ is equal to zero. If at least two coordinates of $\xi$ are equal to zero then the set $\{ (1,0,0), (0,1,0), (0,0,1) (0,0,0)\}$ is a set of representatives in this case. It remains to consider the case when exactly one coordinate of $\xi$ is equal to zero. In view of the cyclic permutation symmetry of the action (in Case 1), it suffices to consider the case when $\xi_3=0$. Then $\l \cdot \xi \in  (  1, K^\times , 0)$ iff $\l \cdot \xi = (1, \xi_1\xi_2\l_1^{-2}, 0)$ (see the case (a)), and the statement (b) follows. 

2(a) $\l \cdot \xi \in (1, K^\times , 1 )$ iff $\frac{\l_3}{\l_1\l_2}\xi_1=1$ and  $\frac{\xi_3}{\l_2\l_3}=1$, and in this case  $\l \cdot \xi =(1,  \frac{\xi_2}{\xi_1\xi_3^2}\l_2^4, 1)$ (since $ \frac{\l_2\xi_2}{\l_1\l_3}= \frac{\xi_2}{\xi_1\xi_3^2}\l_2^4$, use the two equalities); and the statement (a) follows. 

(b) If $\xi \in K^3_{{\rm sing}, 1}(3)$ then $\l \cdot \xi \in (1, K^\times , 0)$ iff $\l \cdot \xi = (1, \xi_1\xi_2\l_1^{-2}, 0)$. Therefore, the set $(1, K^\times / K^{\times 2}, 0)$ is a set of representatives of the orbits in $K^3_{{\rm sing}, 1}(3)$. Clearly, 
 $K^3_{{\rm sing}, 1}(1)=\mT^3 (0,1,1)$  and 
 $K^3_{{\rm sing}, 1}(2)=\mT^3 (1,0,1)$. Now, the statement (b) follows from statement 5.

3(a)   $\l \cdot \xi \in (1,1 , K^\times  )$ iff $\frac{\l_3}{\l_1\l_2}\xi_1=1$ and  $\frac{\xi_2}{\l_1\l_3}=1$, and in this case  $\l \cdot \xi =(1, 1,  \frac{\xi_3}{\xi_1\xi_2^2}\l_1^3)$ (since $ \frac{\xi_3}{\l_2\l_3}= \frac{\xi_3}{\xi_1\xi_2^2}\l_1^3$, use the two equalities); and the statement (a) follows. 

(b) Straightforward.

4(a)  $\l \cdot \xi \in (1,1 , K^\times  )$ iff $\frac{\xi_1}{\l_1\l_2}=1$ and  $\frac{\xi_2}{\l_1\l_3}=1$, and in this case  $\l \cdot \xi =(1, 1,  \frac{\xi_3}{\xi_1\xi_2}\l_1^2)$ (since $ \frac{\xi_3}{\l_2\l_3}= \frac{\xi_3}{\xi_1\xi_2}\l_1^2$, use the two equalities); and the statement (a) follows. 

(b) Straightforward.

5. Statement 5 is obvious.  $\Box $

\begin{corollary}\label{aXA15Mar15}%\marginpar{aXA15Mar15}
Suppose that $\sqrt{K}\subseteq K$ and $ \sqrt[3]{K}\subseteq K$.  Then  for each of the four actions as in Proposition \ref{XA17Mar15}, the map $K^3/\mT^3\ra \{ 0,1\}^3$, $ \mT^3\xi \mapsto \supp (\xi )$ is a bijection with the inverse $\eta \mapsto \mT^3\eta$. 
\end{corollary}

{\bf Proof of Theorem \ref{17Mar15}.}   By writing the equalities (\ref{XBQA1})--(\ref{XBQA3}) in a `cyclicly symmetric way' ($x_2x_1=qx_1x_2+\cdots$, $x_1x_3=qx_3x_1+\cdots$ and $ x_3x_2=qx_2x_3+\cdots $) we see that up to cyclic permutation of the canonical generators there are four cases:\\

{\em Case 1:} $c\neq 0$,  $\beta \neq 0$, $\l \neq 0$. \\

{\em Case 2:} $c\neq 0$,  $\beta \neq 0$, $\l = 0$. \\

{\em Case 3:} $c\neq 0$,  $\beta = 0$, $\l = 0$.\\

{\em Case 4:} $c =  0$, $\beta =  0$, $\l  =  0$.\\

By making the change of variables 
%\marginpar{XpXp}
\begin{equation}\label{XpXp}
x_1'=\l_1^{-1}x_1, \;\; x_2'=\l_2^{-1}x_2, \;\; x_3'=\l_3^{-1}x_3\;\; {\rm where}\;\;  \l_i\in K^\times , 
\end{equation}
the equations (\ref{XBQA1})--(\ref{XBQA3}) can be written as 

%\marginpar{PXBQA1}
\begin{equation} \label{PXBQA1}
x_2'x_1'=qx_1'x_2'+\frac{\l_3}{\l_1 \l_2}cx_3'+\frac{1}{\l_1 \l_2}b_1,
\end{equation} 
%\marginpar{XBQA2}
\begin{equation}\label{PXBQA2}
x_3'x_1'=q^{-1}x_1'x_3'+\frac{\l_2}{\l_1 \l_3}\beta x_2'+\frac{1}{\l_1 \l_3}b_2,
\end{equation}
%\marginpar{PXBQA3}
\begin{equation}\label{PXBQA3}
x_3'x_2'=qx_2'x_3'+\frac{\l_1}{\l_2 \l_3}\lambda x_1'+\frac{1}{\l_2 \l_3}b_3.
\end{equation}

In each of the four cases the following triples \\

Case 1: $(c,\beta , \l )$, \\

Case 2:  $(c,\beta ,b_3 )$, \\

Case 3:  $(c,b_2 ,b_3 )$, \\

Case 4: $(b_1,b_2 ,b_3 )$ \\

are transformed under the change of variables (\ref{XpXp}) as the four $\mT^3$-actions above (see Cases 1--4 of $\mT^3$-actions). Now, using Proposition \ref{XA17Mar15} and using $S_3$-action if necessary we see that (up to $G_3'$) the parameters $c,\beta , \l , b_1, b_2, b_3$ belong precisely to one of the four cases of the theorem (that correspond to Cases 1--4 above) 

If, in addition,  $\sqrt{K}\subseteq K$ and $ \sqrt[3]{K}\subseteq K$ then  the four cases are reduced to ceses 1--4 at the and of the theorem. $\Box$ 

\subsection{Case 2: $q_1-q_3=0$ and $1-q_1q_2\neq 0$}

{\bf Proof of Theorem \ref{QCase-2}.} In view of (\ref{abna}), the equalities (\ref{$(X_1X_1)$})--(\ref{$(X_3X_3)$}) can be written as $\l =0$, $0=0$ and $c=0$, respectively. Then the equalities (\ref{$(X_1)$})--(\ref{$(X_3)$}) can be written as $b_3=0$, $0=0$ and $b_1=0$, respectively. Now, the equality (\ref{$(1)$}) is $0=0$. So, the defining relations (\ref{abna1})--(\ref{abna3}) are as in the theorem. Now, using the change of variables $x_1'=\l_1x_1$, $x_2'=\l_2x_2$ and $x_3'= x_3$ (where $\l_1, \l_2\in K^\times$), we can assume that $\beta , b_2\in \{ 0,1\}$. $\Box$

\subsection{Case 3: $q_1-q_3\neq 0$, $1-q_1q_2= 0$ and  $1-q_2q_3\neq 0$}

Notice that the third condition (i.e., $1-q_2q_3\neq 0$) follows from the first two.\\

{\bf Proof of Theorem \ref{QCase-3}.} In view of (\ref{abna}), the equalities (\ref{$(X_1X_1)$})--(\ref{$(X_3X_3)$}) can be written as $0 =0$, $\beta =0$ and $c=0$, respectively. Then the equalities (\ref{$(X_1)$})--(\ref{$(X_3)$}) can be written as $0=0$, $(q_1-q_3)b_2=0$ and $(1-q_2q_3)b_1=0$, respectively, and so $b_1=b_2=0$. Then the condition (\ref{$(1)$}) holds.  So, the defining relations (\ref{abna1})--(\ref{abna3}) are as in the theorem. Now, using the change of variables $x_1'=\l_1x_1$, $x_2'=\l_2 x_2$ and $x_3'= x_3$ (where $\l_1, \l_2\in K^\times$), we can assume that $\l , b_3\in \{ 0,1\}$. $\Box$

\subsection{Case 4: $q_1-q_3\neq 0$, $1-q_1q_2\neq 0$ and  $1-q_2q_3=0$}

{\bf Proof of Theorem \ref{QCase-4}.} In view of (\ref{abna}), the equalities (\ref{$(X_1X_1)$})--(\ref{$(X_3X_3)$}) can be written as $(1-q_1q_2)\l =0$, $(q_1-q_3)\beta  =0$ and $0=0$, respectively, and so $\l =\beta  =0$.  Then the equalities (\ref{$(X_1)$})--(\ref{$(X_3)$}) can be written as $(q_1q_2-1)b_3=0$, $(q_1-q_3)b_2=0$ and $0=0$, respectively, and so $b_2=b_3=0$.  Then the condition (\ref{$(1)$}) holds.  So, the defining relations (\ref{abna1})--(\ref{abna3}) are as in the theorem. Then the condition (\ref{$(1)$}) holds.   Now, using the change of variables $x_1'=\l_1 x_1$, $x_2'= x_2$ and $x_3'=\l_3x_3$ (where $\l_1, \l_3 \in K^\times$), we can assume that $c , b_1\in \{ 0,1\}$. $\Box$

\subsection{Case 5: $q_1-q_3\neq 0$, $1-q_1q_2\neq 0$ and  $1-q_2q_3\neq 0$}

{\bf Proof of Theorem \ref{QCase-5}.} In view of (\ref{abna}), the equalities (\ref{$(X_1X_1)$})--(\ref{$(X_3X_3)$}) can be written as $\l  =0$, $\beta  =0$ and $c=0$, respectively.  Then the equalities (\ref{$(X_1)$})--(\ref{$(X_3)$}) can be written as $b_3=0$, $b_1=0$ and $b_2=0$, respectively. Then the condition (\ref{$(1)$}) holds. Therefore, $\mA =0$ and $\mB =0$, i.e., $A\simeq \mA^3_Q$.  $\Box$

%%%%%%%%%%%%%%%%%%   Section  7  %%%%%%%%%%%%%%%%%%%%%%%%%%

\section{The bi-quadratic algebras with $q_1=q_2=q_3=1$ (Classification of bi-quadratic algebras $A$ of Lie 
type)}\label{LIETYPE}%\marginpar{LIETYPE}

The aim of this section is to give a classification (up to isomorphism) of the algebras $A=[x_1,x_2,x_3;Q,\mA,\mB]$ of Lie type (Theorem \ref{9Oct18}).\\

{\bf The algebras $A$ of Lie type and their classification.}\\

{\it Definition.} We say that the algebra  $A$ is of {\em Lie type} if $q_1=q_2=q_3=1$.\\

The next theorem  describes all such algebras.
\begin{theorem}\label{8Oct18}%\marginpar{8Oct18}
The algebra $A$ is of Lie type iff $A\simeq U(\CG)/(z-1)$ where $U(\CG)$ is the universal enveloping algebra of a $4$-dimensional Lie algebra $\CG$ such that $z$ is a central nonzero element of $\CG$. 
\end{theorem}
{\it Proof}.\ $(\Leftarrow )$ Let $\CG =Kx_1\bigoplus Kx_2 \bigoplus Kx_3\bigoplus Kz$ be a $4$-dimensional Lie algebra where $\lbrace x_1, x_2, x_3, z\rbrace$ is a $K$-basis of $\CG$.
Let $U=U(\CG)$ be the universal enveloping algebra of  the Lie algebra $\CG$. Then the defining relations of the algebra $U$ are
%\marginpar{zDRel1}
\begin{equation} \label{zDRel1}
x_2x_1-x_1x_2=ax_1+bx_2+cx_3+b_1z,
\end{equation} 
%\marginpar{zDRel2}
\begin{equation}\label{zDRel2}
x_3x_1-x_1x_3=\alpha x_1+\beta x_2+\g x_3+b_2z,
\end{equation}
%\marginpar{zDRel3}
\begin{equation}\label{zDRel3}
x_3x_2-x_2x_3=\lambda x_1+\mu x_2+\nu x_3+b_3z,
\end{equation}
%\marginpar{zDRel4}
\begin{equation}\label{zDRel4}
x_iz-zx_i=0 \  {\rm for} \  i=1,2,3.
\end{equation}
The element $z$ belongs to the centre of the algebra $U$. Then the factor algebra $U/(z-1)$ is the algebra of Lie type with defining constants $a,\ldots, b_3$.

$(\Rightarrow )$ Given an algebra $A$ of Lie type. It can be seen as a Lie algebra $(A,[\cdot,\cdot])$ where $[u,v]:=uv-vu.$ Recall that $A=\bigoplus_{ \alpha\in \N^3} K x^\alpha $ (Theorem \ref{11Apr18}). Hence, the $4$-dimensional vector space $\CG =Kx_1\bigoplus Kx_2 \bigoplus Kx_3\bigoplus Kz$ where $z=1$ is a Lie subalgebra of $A$. Furthermore, there is an algebra epimorphism $U(\CG)\rightarrow A$, $x_i\mapsto x_i$, $z\mapsto 1$. The element $z$ belongs to the center of the algebra $U(\CG)$. By the PBW Theorem and Theorem \ref{11Apr18}, $U/(z-1)\simeq A.$ $\Box $\\

Let $K$ be an algebraically closed field of characteristic zero and $\CG$ be a Lie algebra over $K$. A Lie algebra $\CG$ is called an {\em abelian} Lie algebra if $[\CG , \CG]=0$. The Lie algebra 
%\marginpar{gn2Lie}
\begin{equation*}\label{gn2Lie}
\gn_2=\langle x,y\ |\ [x,y]=y\rangle
\end{equation*} 
is the only (up to isomorphism) $2$-dimensional non-abelian Lie algebra (if $[x,y]=\lambda x+\mu y$ where $\lambda , \mu \in K$ and $\mu \neq 0$, then $[\mu^{-1} x, \lambda x + \mu y]= \lambda x +\mu y;$ it $\mu =0$ then $\lambda\neq 0$ and $[-\lambda^{-1}y,x]=x)$.

Theorem \ref{A12Oct18} is a classification of simple $3$-dimensional Lie algebras. The classification is known (and can be easily proven directly).
\begin{theorem}
\label{A12Oct18} %\marginpar{A12Oct18}
Let $K$ be an algebraically closed field of characteristic zero and $\CG$ be a $3$-dimensional Lie $K$-algebra. Then (up to isomorphism) the Lie algebra $\CG$ is one of the following Lie algebras:
\begin{enumerate}
\item$sl_2 (K)=\langle x,y,h \ | \ [h,x]=2x, \ [h,y]=-2y, \ [x,y]=h\rangle$,
\item $\CH_3= \langle x,y, z\ | \ [x,y]=z, \ [x,z]=0, \ [y,z]=0\rangle$, the $3$-dimensional Heisenberg Lie algebra,
\item $\gn_2 \times Kz$, a dirct product of Lie algebras (where $Kz$ is an abelain $1$-dimensional Lie algebra),
\item an abelian $3$-dimensional Lie algebra.
\end{enumerate}
\end{theorem}

Theorem (\ref{A8Oct18}) is a classification of $4$-dimensional Lie algebras with nontrivial centre.\\

{\bf Classification of $4$-dimensional Lie algebras with nontrivial centre.}

\begin{lemma}
\label{A8Oct18} %\marginpar{A8Oct18}
Let $\CG$ be a Lie algebra with $Z(\CG)\neq0$ and $c \in Z(\CG)\setminus \lbrace0\rbrace$. Then $Kc$ is an ideal of the Lie algebra $\CG$ and there is a short exact sequence of Lie algebras
%\marginpar{KcGG}
\begin{equation*}\label{KcGG}
0\rightarrow Kc\rightarrow \CG \rightarrow \overline{\CG}:={\CG/Kc}   \rightarrow 0. 
\end{equation*}
We can assume that $\CG = \overline{\CG} \oplus Kc,$ a direct sum of vector spaces, i.e., we fix an embedding of the vector space $\overline{\CG}$ into $\CG$, $\overline{\CG} \rightarrow \CG$, $v\mapsto \overline{v}$. Then, for all $u,v \in \CG$, $[\overline{u}, \overline{v}]=[\overline{u,v]}+ \lambda_{u,v}c$ where $\lambda_{\cdot,\cdot}:\overline{\CG}\times \overline{\CG}\rightarrow K, (u,v)\mapsto\lambda_{u,v}$ is a skew-symmetric bilinear form.  
\end{lemma}
{\it Proof.} Straightforward. $\Box$

Theorem \ref{A9Oct18} describes all the $4$-dimensional Lie algebras with nontrivial  centre. We will see that all of them are non-isomorphism Lie algebras, Theorem {\ref{B9Oct18}.
\begin{theorem}
\label{A9Oct18} %\marginpar{A9Oct18}
Let $\CG$ be a $4$-dimensional Lie algebra over a field of characteristic zero such that its centre $Z(\CG)$ is not equal to zero. Fix a nonzero element $c$ of $Z(\CG)$ and let $\overline{\CG} = \CG /Kc$, a Lie factor algebra of dimension $3$.
\begin{enumerate}
\item If \ $\overline{\CG}\simeq sl_2 (K)$ then $\CG \simeq sl_2(K)\times Kc$.
\item If \ $\overline{\CG}\simeq \CH_3$ is the $3$-dimensional Heisenberg Lie algebra, $\CH_3=\langle x,y,z\ |\ [x,y]=z,\ [x,z]=0 ,\ [y,z]=0\rangle$ then either $\CG \simeq \CH_3 \times Kc$ or otherwise $\CG\simeq \CN$ where $\CN=\langle x,y,z,c\ |\ [x,y]=z,\ [x,z]=0,\ [y,z]=c\rangle$ and $c$ is a central element of $\CN$.
\item If \ $\overline{\CG}\simeq \gn_2 \times Kz $ then either $\CG= \gn_2\times Kz\times Kc$ or otherwise $\CG \simeq \CM:= \langle x,y,z,c\ |\ [x,y]=y,\ [x,z]=c,\ [y,z]=0\rangle$ and $c$ is a central element of $\CG$.
\item If \ $\overline{\CG}$ is  an abelian $3$-dimensional Lie algebra then either $\CG$ is an abelian $4$-dimensional Lie algebra or otherwise $\CG \simeq \CH_3 \times Kd$ where $Kd$ is an abelian $1$-dimensional Lie algebra.    
\end{enumerate}   
\end{theorem}
{\it Proof}. $1$. Fix an embedding  of the vector space $sl_2(K)$ into $\CG$, then $$[\overline{h},\overline{x}]=2\overline{x}+2 \lambda_1c,\ [\overline{h},\overline{y}]=-2\overline{y}-2\lambda_2c  ,\ [\overline{x}, \overline{y}]=\overline{h}+\lambda_3c \ {\rm for \  some} \   \lambda_i \in K.$$ Then the linear span of the vectors $x= \overline{x} + \lambda_1 c, \ y=\overline{y}+\lambda_2c$ and $h=h+\lambda_3c$ is a Lie subalgebra of $\CG$ which is isomorphic to $sl_2$ since $[h,x]=2x,\ [h,y]=-2y$ and $[x,y]=h$. We denote this Lie subalgebra by $sl_2(K)$. Then $\CG= sl_2(K)\times Kc$ is a direct product of Lie algebras.

$2$. Fix an embedding of the vector space $\CH_3$ into $\CG$.
Then $$[\overline{x},\overline{y}]=\overline{z}+ \lambda_3c,\ [\overline{x},\overline{z}]=\lambda_1c, \ [\overline{y},\overline{z}]=\lambda_2c \ {\rm for \  some} \   \lambda_i \in  K.$$ Replacing $\overline{z}$ by $\overline{z} +\lambda_3c$, we can assume that $\lambda_3=0$. 

If $\lambda_1=\lambda_2=0$ then the $3$-dimensional vector space $K\overline{x}\oplus K\overline{y} \oplus K\overline{z}$ is a Lie subalgebra of $\CG$ which is isomorphic to $\CH_3$, and $\CG \simeq \CH_3 \times Kc$ is a direct product of Lie algebras.

 If $\lambda_2\neq0$ then replacing $c$ by $\lambda_2c$ we can assume that $\lambda_2=1$. Then replacing $\overline{x}$ by $ \overline{x}-\lambda_1\overline{y}$ we can assume that $\lambda_1=0$, i.e., $\CG\simeq \CN$. 
 
 If $\lambda_2=0$ and $\lambda_1 \neq 0$ then  the linear transformation $\overline{x}\mapsto\overline{y},\ \overline{y} \mapsto -\overline{x}, \ \overline{z}\mapsto \overline{z},\  c\mapsto c $ we come to the previous case (where $\lambda_2\neq0$).

$3.$ Fix an embedding of the vector space $\gn_2\times Kz$ into $\CG$. Then $$[\overline{x},\overline{y}]=\overline{y}+\lambda_3c, \ [\overline{x},\overline{z}]=\lambda_1c, \ [\overline{y},\overline{z}]=\lambda_2c \  {\rm \ for \   some} \  \lambda_i \in K.$$ Replacing $ \overline{y}$ by $\overline{y}+\lambda_3c$ we can assume the $\lambda_3=0$.

If $\lambda_1=\lambda_2=0$ then the $3$-dimensional vector space $K \overline{x}\oplus K \overline{y}\oplus K \overline{z}$ is a Lie subalgebra of $\CG$ which is isomorphic to $\gn_2 \times Kz$, and so $\CG \simeq \gn_2 \times Kz\times Kc$ is a direct product of Lie algebras. 

If, say, $\lambda_1\neq0$ then replacing $c$ by $\lambda_1 c$ we can assume the $\lambda_1=1$. Then replacing $\overline{y}$ by $\overline{y} - \lambda_2\bx$ we can assume these $\lambda_2=0$, i.e., $\CG\simeq \CM$.

 If $\lambda_1=0$ and $\lambda_2 \neq 0$ then replacing $c$ by $\lambda_2c$ we can assume the $\lambda_2=1$, i.e., $[\overline{x},\overline{y}]=\overline{y}, \ [\overline{x},\overline{z}]=0, \  [\overline{y},\overline{z}]=c$. This is not a Lie bracket since
$0=[\by,[\bx,\bz]]=[[\by,\bx],\bz]+[\bx,[\by,\bz]]$
$=-[\by,\bz]+[\bx,c]=-c$, a contradiction.

$4.$ Fix an embedding of the $3$-dimensional abelian Lie algebra $\overline{\CG}=\langle\bx_1,\bx_2,\bx_2\rangle$ into $\CG$ as vector space. We assume that the Lie algebra $\CG$ is not abelian. To finish the proof of statement $4$ we have to show that $\CG\simeq \CH_3 \times Kd$.
Notice that $$[\bx_1,\bx_2]=\lambda c, \ [\bx_1,\bx_3]=\mu c, \ [\bx_2, \bx_3]=\delta c.$$ Without less of generality we may assume that $\lambda \neq 0$ since $\CG$ is not an abelian Lie algebra. Then replacing $\bx_1$ by $\lambda^{-1} \bx_1$ we can assume that $\lambda =1$. Hence the vector space $$\langle \bx_1 , \bx_2,c\  |\  [\bx_1,\bx_2]=c\rangle$$ is a $3$-dimensional Heisenberg Lie algebra $\CH_3$. Then replacing  $\bx_3$  by $\bx_3-\mu \bx_2$ we can assume that $\mu=0$, i.e., $$[\bx_1,\bx_2]=c , \ [\bx_1,\bx_3]=0 \ {\rm\  and}  \  [\bx_2,\bx_3]=\d c.$$ Then either $\d=0$ and in this case $\CG \simeq \CH_3 \times Kx_3$ is a direct product of Lie algebras or $\d \neq 0$. If $\d\neq 0$  then changing $\bx_3$ to $\d^{-1} \bx_3$ we can assume that $\d=1$, i.e., $[\bx_2,\bx_3]=c$. Then $[\bx_2,\bx_3+ \bx_1]=c-c=0.$ Then changing $\bx_3$  to $\bx_3 + \bx_1$  we can assume the $\d =0$, and we are in the situation  of the previous case. $\Box$

The next theorem is a classification of $4$-dimensional Lie algebras with nonzero centre.
\begin{theorem}
\label{B9Oct18}%\marginpar{B9Oct18}
Let $\CG$ be a $4$-dimensional Lie algebras over an algebraically closed field of characteristic zero such that $Z(\CG)\neq0$. Fix a nonzero element $c$ of $Z(\CG)$. Then $\CG$ is isomorphic to one of  the following non-isomorphic Lie algebras:
\begin{enumerate}
\item $\CG$ is an abelian $4$-dimensional Lie algebra,
\item $sl_2(K)\times Kc$,
\item $\CH_3\times Kc$,
\item $\CN= K\langle x,y,z,c\ |\ [x,y]=z, \  [x,z]=0,\  [y,z]=c\rangle$ and $c$ is a central element of the Lie algebra $\CN$.
\item $\gn_2 \times Kz \times Kc$.
\item $\CM=\langle x,y,z,c \ |\ [x,y]=y,\  [x,z]=c,\  [y,z]=0 \rangle $ and $c$ is a central element of the Lie algebra $\CM$.
 
\end{enumerate} 
\end{theorem}

{\it Proof}. It follows from Theorem \ref{A9Oct18} and the classification of $3$-dimensional Lie algebras (Theorem \ref{A12Oct18}) that the Lie algebra $\CG$ is isomorphic to one of the algebras in statements 1-6. We denote by $\CG_i$ the Lie algebras in the statement $i$. Let $n=\dim\  Z(\CG)$.
Then
\begin{eqnarray*}
 n=1&:& \CG_2,\  \CG_4,\ \CG_6,\\
n=2&:&\CG_3,\ \CG_5,\\
 n=4&:& \CG_1.
\end{eqnarray*}
So, it remains to show that the Lie algebras in the three cases above are not isomorphic. The Lie algebras $\CG_3$ and $\CG_4$ are nilpotent. The Lie algebras $\CG_5$ and $\CG_6$ are solvable but non-nilpotent. The Lie algebra $\CG_2$ is a direct product of a simple Lie algebra  and an abelian Lie algebra. Hence, all Lie algebras $\CG_1, \ldots , \CG_6$ are not isomorphic. $\Box$\\

{\bf Classification (up to isomorphism) of algebras of Lie type.} Theorem \ref{9Oct18} is such a classification.\\

{\bf Proof of  Theorem \ref{9Oct18}}. It follows from Theorem \ref{8Oct18} and Theorem \ref{B9Oct18} that an algebra $A$ of Lie type  is  isomorphic to one of the algebras in statements 1-6.

 Let $\CA_i$ be the algebra in the statement $i$ for $i=1,\ldots,6$. It remains to prove that the algebras $\CA_1 - \CA_6$ are not isomorphic. Since the algebra $\CA_1$ is the only commutative algebra out of six, we have to show that the algebras $\CA_2-\CA_6$ are not isomorphic. In Proposition \ref{A13Oct18}, we collect isomorphism-invariant properties
of the algebras $\CA_2-\CA_6$ that show that they are not isomorphic.   By Proposition  \ref{A13Oct18}, the algebra $\CA_6$ is the only algebra out of $\CA_2, \ldots, \CA_6$ that has trivial centre. Hence, it remains to show that the algebra $\CA_2-\CA_5$ are not isomorphic. By Proposition \ref{A13Oct18}, the algebra $\CA_2$ is the only algebra out of $\CA_2 - \CA_5$ that has simple $2$-dimensional module. So, it remains to show that the algebras $\CA_3,\CA_4$ and $\CA_5$ are not isomorphic. By Proposition \ref{A13Oct18}, the centre of each of the algebras $\CA= \CA_3, \CA_4, \CA_5$ is a polynomial algebra $K[t] $ in a single a variable. The algebra $\CA=\CA_4$ is the only algebra of $\CA_3,\CA_4$ and $\CA_5$ that has the property  that all factor algebras $\CA/(t-\lambda)$, where $\lambda \in K$, are simple (since the Weyl algebra $A_1$ is a simple algebra). So, it remains to show that $\CA_3 \not\simeq\CA_5$. This follows from statements $2(c)$ and $4(c)$ of Proposition \ref{A13Oct18}. $\Box$  
\begin{proposition}\label{A13Oct18}%\marginpar{A13Oct18}
Let $K$ be  an algebraically closed field of characteristic zero
\begin{enumerate}
\item Let $\CA_2=U(sl_2(K)).$ Then 
\begin{enumerate}
\item $Z(\CA_2)=K[C]$ where $C$ is the Casimir element.
\item For each national number $n=1,2,\ldots$ there is a unique (up to isomorphism) simple $n$-dimensional $\CA_2$-module.
\end{enumerate}
\item Let $\CA_3=U(\CH_3)$. Then
\begin{enumerate}
\item $Z(\CA_3)=K[z].$
\item $\CA_3/(z)\simeq K[x,y]$ is a polynomial algebra and $\CA_3(z-\lambda)\simeq A_1$ is the Weyl algebra for all $\lambda \in K^\times$.
\item $\{K_{\mu,\d}:=\CA_3/(x-\mu,y-\d, z)\ |\  \mu,\d\in K\}$ is the set of isomorphism classes of simple finite dimensional $\CA_3$-modules and $\dim (K_{\mu,\d})=1$ and  $\ann_{Z(\CA_3)}(K_{\mu,\d})=zK[z]$ for all $\mu, \d \in K$. 
\end{enumerate}
\item Let $\CA_4= U(\CN)/(c-1)=A_1 \otimes K[x'].$ Then 
\begin{enumerate}
\item $Z(\CA_n)=K[x']$ where $x'=x+\dfrac{1}{2}z^2$.
\item All simple $\CA_4$-modules are infinite dimensional.
\item For all $\lambda \in K, \CA_4 /(x'-\lambda)\simeq A_1$, the Weyl algebra.
\end{enumerate}
\item Let $\CA_5 =U(\gn_2\times Kz)=U(\gn_2)\otimes K[z]$. Then
\begin{enumerate}
\item $Z(\CA_5)=K[z]$.
\item For all $\lambda \in K, \CA_5 /(z-\lambda)\simeq U(\gn_2)\simeq K[x][y;\sigma]$ is a skew polynomial algebra where $\sigma(x)=x+1$.
\item $\{K_\mu(\lambda):=\CA_5/(x-\mu,y,z-\lambda)\ |\  \lambda,\mu \in K\}$ is the set of isomorphism classes of simple finite dimensional $\CA_5$-modules and $\dim (K_\mu(\lambda))=1$ for all $\mu,\lambda \in K$, and $\ann_{Z(\CA_5)}(K_\mu(\lambda))= (z-\lambda)K[z]$ for all $\lambda,\mu \in K$.   
\end{enumerate}
\item Let $\CA_6:=U(\CM)/(c-1)\simeq A_1[y,\sigma]$ where $A_1=K\langle x,z|[x,z]=1\rangle$ is the Weyl algebra, $\sigma(x)=x+1$ and $\sigma(z)=z$. Then 
\begin{enumerate}
\item $Z(\CA_6)=K.$
\item There is no simple finite dimensional $\CA_6$-modules.
\end{enumerate}  
\end{enumerate}
\end{proposition}
{\it Proof.} $1$. Statement $1$ is a well-known fact.

$2$.(a) The statement (a) is a well-known (and easy to prove) fact.

(b) The statement (b) is obvious.

(c) The statement (c) follows from the statement (b) and the fact then the Weyl algebra $A_1$ is a simple infinite dimensional algebra (and therefore, all simile $A_1$-modules are infinite dimensional).

$3.$ (a) $Z(\CA_4)=Z(A_1\otimes K[x'])=Z(A_1)\otimes Z(K[x'])=K\otimes K[x']=K[x']$ since $Z(A_1)=K$.

(b),(c): The statement (c) is obvious. The statement (b) follows from the statement (c).

$4.$ (a) $Z(\CA_5 )= Z(U(\gn_2 \times Kz))= Z(U(\gn_2) \otimes K[z])= Z(U(\gn_2))\otimes K[z]=K[z].$

(b) The statement (b) follows from the statement (a).

(c) The statement (c) follows from the statement (b). In more detail, the algebra $U(\gn_2)\simeq K[x][y;\sigma]$ is a skew polynomial algebra. The element $y$ is a regular normal element $(yU(\gn_2)=U(\gn_2)y)$. The localization of the algebra $U(\gn_2)$ at the powers of the element $y$ is a simple skew polynomial ring $B_1=K[x][y,y^{-1};\sigma]$. If $M$ is a simple $U(\gn_2)$-module the either $yM=0$ or otherwise the map $y: M\rightarrow M, m\mapsto ym$ is a bijection and therefore, the $ U(\gn_2)$-module, $M$ is also a simple $B$-module, and so $\dim (M)=\infty$. So, if $M$ is a simple finite dimensional $U(\gn_2)$-module then $yM=0$, and so $M$ is a simple finite dimensional module over the polynomial algebra $U(\gn_2)/(y) \simeq K[x]$. Now, the statement (b) is obvious. 

$5.$ (a) The centralize $C(y)$ of the element $y$ in $\CA_6$ is equal to the polynomial algebra $K[y,z]$ (since $\sigma(x)=x+1$ and $\sigma(z)=z$). Similarly, the centralize $C(x)$ of $x$ is the polynomial algebra $K[x]$. Since $Z(\CA_6) \subseteq  C(x)\cap C(y)=K[x]\cap K[y,z]=K,$ we must have $ Z(\CA_6)=K$. 

(b) Since the Weyl algebra $A_1$ is a subalgebra of $\CA_6$, every nonzero $\CA_6$-module is also a nonzero $A_1$-module. The Weyl algebra $A_1$  is a simple infinite dimensional algebra. So, all nonzero $A_1$-modules are infinite dimensional, and the statement (b) follows.  $\Box$

\small{

Department of Pure Mathematics

University of Sheffield

Hicks Building

Sheffield S3 7RH

UK

email: v.bavula@sheffield.ac.uk
}

\end{document}